\documentclass{article}

\usepackage{vmargin,epsfig}
\usepackage[english,francais]{babel}
\usepackage[latin1]{inputenc}
\usepackage[tbtags]{amsmath}
\usepackage{amsthm,amssymb}
\usepackage{pifont}
\usepackage[all]{xy}

\newtheorem{theo}{Théorème}[section]
\newtheorem{lemme}[theo]{Lemme}
\newtheorem{prop}[theo]{Proposition}
\newtheorem{cor}[theo]{Corollaire}
\theoremstyle{definition}
\newtheorem{deftn}[theo]{Définition}
\newtheorem{question}[theo]{Question}

\def\leq{\leqslant}
\def\geq{\geqslant}

\def\N{\mathbb{N}}
\def\Z{\mathbb{Z}}
\def\Q{\mathbb{Q}}

\def\F{\mathbb{F}}
\def\O{\mathcal{O}}

\def\pa#1{\left(#1\right)}

\def\epsilon{\varepsilon}

\def\calD{\mathcal{D}}
\def\calL{\mathcal{L}}
\def\calM{\mathcal{M}}
\def\calN{\mathcal{N}}

\def\id{\text{\rm id}}

\def\Fil{\text{\rm Fil}}
\def\Gal{\text{\rm Gal}}
\def\pgcd{\text{\sc Pgcd}}

\def\ModphiN#1{\text{\rm Mod}^{\phi,N}_{/#1}}
\def\Modphi#1{\text{\rm Mod}^{\phi}_{/#1}}
\def\MF{\text{\rm MF}}
\def\calMF{\mathcal{MF}}
\def\fa{\text{\rm fa}}

\def\st{\text{\rm st}}

\def\tr{\text{\rm tronc}}
\def\H{\mathbf{H}}

\def\barM{\overline{\mathcal{M}}}
\def\calA{\mathcal{A}}
\def\e{e}

\title{Polygones de Hodge, de Newton et de l'inertie modérée\\
des représentations semi-stables}
\author{Xavier Caruso et David Savitt}
\date{Juin 2008}

\begin{document}

\maketitle

\setcounter{tocdepth}{2}
\tableofcontents

\bigskip

\noindent \hrulefill

\bigskip

\section{Introduction}

Soient $p$ un nombre premier, et $k$ un corps parfait de caractéristique
$p$. On note $W = W(k)$ l'anneau des vecteurs de Witt à coefficients
dans $k$, $K_0$ le corps des fractions de $W$ et on fixe $K$ une
extension totalement ramifiée de $K_0$ de degré $e$, d'anneau des
entiers noté $\O_K$. On fixe également $\bar K$ une cloture algébrique
de $K$, on désigne par $G_K$ le groupe de Galois absolu de $K$ et par
$I$ le sous-groupe d'inertie.

\medskip

Si $V$ est une $\Q_p$-représentation semi-stable (voir \cite{fontaine})
de $G_K$ de dimension $d$, on sait lui associer plusieurs invariants
numériques. Il y a classiquement les poids de Hodge-Tate et les pentes
de l'action du Frobenius sur le module filtré de Fontaine correspondant.
Ce sont tous deux des $d$-uplets de rationnels (entiers pour les poids
de Hodge-Tate) que l'on représente parfois sous la forme de polygones
appelés alors respectivement polygone de Hodge et polygone de Newton.
(Pour les définitions de ces multiples objets, voir paragraphe
\ref{subsec:fontaine}.) On sait que ces deux polygones ont même point
d'arrivée et que le polygone de Hodge est toujours situé en-dessous
de celui de Newton. Si $h_1 \leq h_2 \leq \cdots \leq h_d$ sont les poids
de Hodge-Tate et $n_1 \leq n_2 \leq \cdots \leq n_d$ les pentes de
Frobenius, la condition visuelle sur les polygones se traduit par les
inégalités
$$h_1 + \cdots + h_k \leq n_1 + \cdots + n_k$$
pour tout $k \in \{1, \ldots, d\}$ avec égalité si $k=d$.

\medskip

Il existe toutefois un troisième invariant numérique que l'on peut
associer à $V$. Rappelons avant de le définir que toute
$\F_p$-représentation irréductible $H$ du groupe d'inertie $I$ est
simplement décrite par ses poids de l'inertie modérée : ils forment une
suite de $d' = \dim H$ entiers compris entre $0$ et $p-1$ (définie à
permutation circulaire près). (Pour des précisions sur cette
classification, le lecteur pourra se reporter au paragraphe 1 de
\cite{serre}.)

Considérons à présent $T$ un $\Z_p$-réseau dans $V$ stable par l'action
du groupe de Galois (un tel réseau existe toujours par compacité de
$G_K$), et focalisons-nous sur la représentation $T/pT$ restreinte au
groupe d'inertie $I$. D'après les rappels que nous venons de faire, à
tout quotient de Jordan-Hölder de cette représentation, il correspond
une suite d'entiers compris entre $0$ et $p-1$. D'autre part, un
théorème de Brauer-Nesbitt (voir paragraphe 82.1 de \cite{curtis})
montre que la famille des quotients de Jordan-Hölder de $T/pT$ ne dépend
que de la représentation $V$ (et pas du choix de $T$).  Ainsi, en
regroupant les poids associés à chacun des quotients de Jordan-Hölder,
on obtient une collection de $d$ entiers compris entre $0$ et $p-1$ qui
est \emph{canoniquement} associée à $V$. Notons-les $i_1 \leq i_2 \leq
\cdots \leq i_d$.

Il est naturel de se demander si ces entiers ont un rapport avec les
autres invariants.  Lorsque $e=1$ et lorsque la représentation est
cristalline, Fontaine et Laffaille ont montré dans
\cite{fontaine-laffaille} que les poids de l'inertie modérée sont les
mêmes que les poids de Hodge-Tate. Toutefois, ce résultat simple est mis
en défaut pour les représentations semi-stables non cristallines, même
dans les cas les plus simples : dimension $2$ et $e = 1$. Plus
précisément, dans cette situation, Breuil et Mézard ont montré dans
\cite{breuil-mezard} par un calcul explicite que si les poids de
Hodge-Tate sont $0$ et $r$, alors les poids de l'inertie modérée peuvent
être n'importe quel couple $(i_1, i_2)$ pourvu que $i_1 + i_2 = r$.

\medskip

Dans cette note, nous dégageons des contraintes explicites sur les
poids de l'inertie modérée dans le cas général. Commençons pour cela par
définir \emph{le polygone de l'inertie modérée} comme le polygone
associé aux rationnels $i'_k = \frac {i_k} e$. Nous montrons le théorème
suivant:

\begin{theo}
\label{theo:hodge}
Supposons que tous les poids de Hodge-Tate soient compris entre $0$ et
$r$ pour un entier $r$ vérifiant $er < p-1$. Alors, le polygone de
l'inertie modérée est situé au-dessus du polygone de Hodge. De plus ils
ont même point terminal. Autrement dit, pour tout $k \in \{1, \ldots,
d\}$ : 
$$e(h_1 + \cdots + h_k) \leq i_1 + \ldots + i_k$$
avec égalité si $k = d$.
\end{theo}

\noindent
{\it Remarque.} L'hypothèse $er < p-1$ se justifie aisément puisque les
poids de l'inertie modérée sont par définition bornés par $p-1$. On
pourrait toutefois se demander si le résultat demeure avec la condition
plus faible $er \leq p-1$. La réponse est négative en général, mais il
est tout à fait probable qu'elle devienne positive sous certaines
hypothèses supplémentaires simples.

\bigskip

Forts de ce résultat, nous essayons ensuite de comprendre s'il y a un
rapport entre polygone de l'inertie modérée et polygone de Newton. Nous
obtenons le résultat partiel suivant:

\begin{theo}
\label{theo:newton}
Supposons que tous les poids de Hodge-Tate soient compris entre $0$ et
$r$ pour un entier $r$ vérifiant $er < p-1$. Alors, tant que les
polygones de Hodge et de Newton ne se séparent pas, le polygone de
l'inertie modérée coïncident aussi avec eux. En particulier, si polygone
de Hodge et polygone de Newton coïncident, ils coïncident également avec
le polygone de l'inertie modérée.

Avec des égalités, cela se traduit de la façon suivante. Soit $k \in
\{1, \ldots, d\}$. Supposons que pour tout $k' \leq k$, on ait $h_{k'} =
n_{k'}$. Alors $i_k = e h_k ( = e n_k)$.
\end{theo}

L'organisation de cette note se fait comme suit. La partie
\ref{sec:prelim} est consacrée à des rappels sur la théorie de Breuil
développée dans \cite{breuil-griff}, \cite{breuil-ens},
\cite{breuil-invent} et \cite{caruso-crelle}, qui joue un rôle central
dans la démonstration des résultats précédemment annoncés. La partie
\ref{sec:position}, quant à elle, regroupe les preuves desdits
résultats. Finalement, dans une dernière partie, nous menons à terme le
calcul des invariants précédents sur un exemple simple en dimension $2$.
Il en ressort que, contrairement à ce qui se passe pour $e=1$, les
polygones de Hodge et de l'inertie modérée ne sont en général pas
confondus pour les représentations cristallines.

\paragraph{Remerciements}

Pendant l'accomplissement de ce travail, le second auteur a été
partiellement financé par le \emph{NSF grant DMS-0600871} et lui en est
reconnaissant. De plus, une partie de cet article a été écrite alors que
les auteurs étaient respectivement en visite à l'Université de Bonn et
au \emph{Max-Planck-Institut für Mathematik}. Que ces instituts soient
ici remerciés comme il se doit pour leur accueil chaleureux. C'est
finalement un plaisir pour le premier auteur de remercier Christophe
Breuil de lui avoir conseillé de se pencher sur la problématique de ce
papier, ainsi que pour ses encouragements après les premiers résultats.

\section{La théorie de Breuil}
\label{sec:prelim}

Cette section est dédiée aux préliminaires de théorie de Hodge
$p$-adique, rationnelle et entière. Nous rappelons dans un premier temps
la correspondance, dûe à Fontaine, entre représentations $p$-adiques
semi-stables et $(\varphi, N)$-modules filtrés faiblement admissibles.
Ensuite, nous détaillons la théorie de Breuil qui permet de comprendre
les réseaux à l'intérieur de telles représentations galoisiennes et les
réductions modulo $p$ de ces réseaux.

\subsection{Théorie de Hodge $p$-adique rationnelle}
\label{subsec:fontaine}

\subsubsection*{Les $(\varphi,N)$-modules filtrés de Fontaine}

\paragraph{Définition}

On note $\MF_K(\varphi,N)$ la catégorie dont les objets sont la 
donnée de :
\begin{itemize}
\item un $K_0$-espace vectoriel $D$ de dimension finie ;
\item une filtration décroissante $\Fil^t D_K$ exhautive et séparée sur 
$D_K = D \otimes_{K_0} K$ par des sous $K$-espaces vectoriels ;
\item une application $K_0$-semi-linéaire (par rapport au Frobenius sur $K_0$)
injective $\varphi : D \to D$ (appelée Frobenius) ;
\item une application $K_0$-linéaire $N : D \to D$ (appelée monodromie) vérifiant
$N \varphi = p \varphi N$.
\end{itemize}
Les morphismes dans cette catégorie sont les applications
$K_0$-linéaires compatibles au Frobenius, à la monodromie et à la
filtration après extension des scalaires à $K$. On associe à tout objet
$D \in \MF_K(\varphi,N)$ deux invariants numériques qui sont son nombre
de Hodge et son nombre de Newton. Si $D$ est de dimension $1$, le nombre
de Hodge $t_H(D)$ est défini comme le plus grand entier $t$ tel que
$\Fil^t D_K = D_K$ alors que le nombre de Newton $t_N(D)$ est la pente
du Frobenius sur $D$ (c'est-à-dire la valuation $p$-adique de $\alpha
\in K$ tel que $\phi(x) = \alpha x$ pour un $x \in D$). Si $D$ est de
dimension $d$, on pose par définition $t_N(D) = t_N(\Lambda^d D)$ et
$t_H(D) = t_H(\Lambda^d D)$. Un objet $D$ est dit faiblement admissible
si $t_N(D) = t_H(D)$ et si pour tout $D' \subset D$ stable par $\varphi$
et $N$ et muni de la filtration induite, on a $t_H(D') \leq t_N(D')$. On
note $\MF^\fa_K(\varphi,N)$ la sous-catégorie de $\MF_K(\varphi,N)$
formée des objets faiblement admissibles. On montre que c'est une
catégorie abélienne, stable par produit tensoriel et dualité.

\paragraph{Polygones de Hodge et de Newton}

Plutôt que les nombres de Hodge et de Newton, on considère parfois les
polygones de Hodge et Newton. Comme ceux-ci sont centraux dans toute 
cette note, nous en rappelons la définition.
Rappelons dans un premier temps que si $n_1 \leq n_2 \leq \cdots \leq
n_d$ sont des nombres, on leur associe un \og polygone \fg\ de la façon 
suivante :
on rejoint dans le plan les points de coordonnées $(k, n_1 + \cdots +
n_k)$ et on trace les verticales issues du début et de la fin de la ligne
brisée précédente. Le point $(d, n_1 + \cdots + n_d)$ est appelé
\emph{point d'arrivée} ou \emph{point terminal} du polygone. Si les
entiers $n_i$ ne sont pas triés, on commence par le faire.

Si maintenant $D$ est un objet de $\MF_K(\varphi,N)$, son polygone de
Hodge est le polygone associé aux entiers $t$ pour lesquels $\Fil^{t+1}
D_K \neq \Fil^t D_K$ avec la convention que l'entier $t$ est repété
autant de fois que la dimension de $\Fil^{t+1}D_K / \Fil^t D_K$. Le
polygone de Newton de $D$ quant à lui est le polygone associé aux pentes
de l'action de Frobenius ; c'est aussi le polygone de Newton (usuel) du
polynôme caractéristique de la matrice de $\varphi$ dans une base
quelconque de $D$ sur $K_0$. On montre que si $D$ est faiblement
admissible alors le polygone de Hodge est en-dessous le polygone de
Newton et que ceux-ci ont même point d'arrivée.

\paragraph{Lien avec les représentations semi-stables}

Dans \cite{fontaine2}, Fontaine construit un anneau $B_\st$ muni de
structures supplémentaires (dont nous ne rappelons pas la définition
ici, car elle ne nous servira pas). Il montre que si $V$ est une
représentation semi-stable de $G_K$, alors :
$$D_\st(V) = (B_\st \otimes_{\Q_p} V)^{G_K}$$
est un objet de $\MF^\fa_K(\varphi, N)$. Mieux, le foncteur $D_\st$
établit une anti-équivalence de catégories entre la catégorie des
représentations galoisiennes semi-stables et $\MF^\fa_K(\varphi,N)$. Les
polygones de Hodge et de Newton de la représentation galoisienne sont
alors définis comme les polygones de Hodge et de Newton de l'objet de
$\MF^\fa_K(\varphi,N)$ correspondant.

\subsubsection*{Les $S_{K_0}$-modules filtrés de Breuil}

\paragraph{Définitions}

Soit $f_\pi : W[u] \to \O_K$ l'application qui envoie $u$ sur $\pi$.
Notons $E(u)$ le polynôme minimal de $\pi$ sur $K_0$ de sorte que le
noyau de $f_\pi$ soit l'idéal principal engendré par $E(u)$. On note $S$
le complété $p$-adique de l'enveloppe à puissances divisées de $W[u]$
relativement au noyau de $f_\pi$ (et compatible avec les puissances
divisées canoniques sur $pW[u]$). L'application $f_\pi$ se prolonge à
$S$. De plus, $S$ est naturellement muni d'une filtration $\Fil^i S$ (la
filtration à puissances divisées), d'un Frobenius $W$-semi-linéaire
$\phi : S \to S$ continu pour la topologie $p$-adique défini par
$\phi(u) = u^p$, et d'un opérateur de monodromie $W$-linéaire $N : S \to
S$ continu pour la topologie $p$-adique, vérifiant la condition de
Leibniz et défini par $N(u) = -u$. On note $c = \phi(E(u))/p$ ; c'est
une unité de $S$. Finalement, on désigne par $f_0 : S \to W$ le
morphisme d'anneau (continu pour la topologie $p$-adique) qui envoie $u$
et toutes ses puissances divisées sur $0$.

On pose $S_{K_0} = S \otimes_W {K_0}$ ; les structures supplémentaires que l'on 
vient de définir s'y prolongent.

\medskip

Ceci permet de définir la catégorie $\calMF(\phi, N)$. Ses objets sont
la donnée de :
\begin{itemize}
\item un $S_{K_0}$-module $\calD$ libre de rang fini ;
\item une filtration $\Fil^t \calD$ telle que $\Fil^i S \cdot \Fil^t
\calD \subset \Fil^{t+i} \calD$ pour tous $i$ et $t$ ;
\item un Frobenius $S_{K_0}$-semi-linéaire $\phi : \calD \to \calD$ ;
\item un opérateur de monodromie $N : \calD \to \calD$ vérifiant $N \phi
= p \phi N$ et $N (\Fil^t \calD) \subset \Fil^{t-1} \calD$ pour tout
$t$.
\end{itemize}

\paragraph{Une équivalence de catégories}

Si $D$ est un objet de $\MF_K(\varphi,N)$, on peut lui associer un objet
de $\calD \in \calMF(\phi, N)$ de la façon suivante. On pose $\calD = D
\otimes_{K_0} S_{K_0}$. Il est muni de $\phi = \phi \otimes \varphi$ et
de $N = N \otimes 1 + 1 \otimes N$. La filtration, quant à elle, est
définie par :
$$\Fil^t \calD = \{x \in \calD \, / \, \forall i \geq 0, \, f_\pi
(N^i(x)) \in \Fil^{t-i} D \}.$$

On montre (voir \cite{breuil-griff}, section 6) que cette association
définit une équivalence de catégories entre $\MF_K(\varphi, N)$ et
$\calMF(\phi,N)$. Un quasi-inverse est décrit comme suit. Soit $\calD
\in \calMF(\phi,N)$. On pose $D = \calD \otimes_{S_{K_0}} K_0$ (où le
morphisme $S_{K_0} \to K_0$ est $f_0$). Il est muni des opérateurs
$\phi$ et $N$ qui passent au quotient et définissent respectivement le
Frobenius et l'opérateur de monodromie de l'objet de $\MF_K(\varphi,N)$.
La définition de la filtration demande un peu plus de travail, et
notamment un lemme préliminaire.

\begin{lemme}
Avec les notations précédentes, il existe une unique section $s : D \to
\calD$ de $\id \otimes f_0$, qui soit $K_0$-linéaire et compatible au
Frobenius.
\end{lemme}

\begin{proof}
Voir proposition 6.2.1.1 de \cite{breuil-griff}.
\end{proof}

\noindent
Le section $s$ du lemme fournit en particulier une flèche $K_0$-linéaire
$D \to \calD/\Fil^1 \calD$. En utilisant que $\calD/\Fil^1 S_{K_0}
\calD$ est un espace vectoriel sur $S_{K_0} / \Fil^1 S_{K_0} \simeq K$
(par la flèche $u \mapsto \pi$) dont le rang est le même que celui de
$D$, on montre facilement que $s$ s'étend en un \emph{isomophisme}
$K$-linéaire $s_K : D_K \to \calD/\Fil^1 S_{K_0} \calD$. La filtration
est alors obtenue par $\Fil^t D = s_K^{-1} (\Fil^t \calD / Fil^1 S_{K_0}
\calD)$.

Finalement, mentionnons que de la définition de la filtration, on déduit
l'égalité suivante, qui nous sera utile pour la suite :
\begin{equation}
\label{eq:filt}
\Fil^{t-i} \calD = \{x \in \calD \, / \, E(u)^i x \in \Fil^t \calD \}
\end{equation}
valable pour tous entiers $t$ et $i$.

\subsection{Théorie de Hodge $p$-adique entière et de torsion}
\label{subsec:torsion}

À partir de maintenant, on se donne un entier positif $r < p-1$ et on se
restreint aux représentations pour lesquelles $\Fil^0 D_K = D_K$ et
$\Fil^{r+1} D_K = 0$ sur le module filtré de Fontaine $D$ associé.
L'intérêt d'utiliser la théorie de Breuil est qu'elle permet de décrire
les réseaux dans les représentations semi-stables (vérifiant l'hypothèse
précédente).

\paragraph{Modules fortement divisibles}

Soit $\calD$ un objet de $\calMF(\phi,N)$ associé à un objet $D \in
\MF_K (\varphi,N)$ dont la filtration est comprise entre $0$ et $r$. Un
\emph{module fortement divisible} (ou \emph{réseau fortement divisible})
dans $\calD$ est un sous-$S$-module $\calM \subset \calD$ vérifiant :
\begin{itemize}
\item $\calM$ est libre de rang fini sur $S$ ;
\item le morphisme naturel $\calM \otimes_S S_{K_0} \to \calD$ est un
isomorphisme ;
\item $\phi(\Fil^r \calM) \subset p^r \calM$ où $\Fil^r \calM = \calM
\cap \Fil^r \calD$ ;
\item $\frac{\phi}{p^r}(\Fil^r \calM)$ engendre $\calM$ sur $S$.
\end{itemize}
On montre (voir corollaire 2.1.4 de \cite{breuil-invent}) que si $\calD$
admet un réseau fortement divisible, alors $\calD$ correspond à un $D$
faiblement admissible $V$, et donc à une représentation galoisienne. La
réciproque est également vraie : si $er < p-1$ (cas qui nous intéresse
principalement ici), c'est le résultat principal de
\cite{breuil-invent}, sinon c'est une conséquence des travaux de Kisin
(voir \cite{kisin05}) et de Liu (voir \cite{liu}). En outre, dans le cas
où $D$ est faiblement admissible, la dernière condition (à savoir \og
$\frac{\phi}{p^r}(\Fil^r \calM)$ engendre $\calM$ \fg) est impliquée par
les précédentes, ce qui explique par exemple que nous ne la vérifierons
dans la section \ref{sec:exemple}.

D'autre part, on dispose d'un foncteur $T_\st$ (dont nous renvoyons à
\cite{breuil-ens} pour la définition) qui à $\calM$ associe un
$\Z_p$-réseau stable par Galois $T$ à l'intérieur de $V$. Ce foncteur
établit en fait une bijection entre les réseaux fortement divisibles
dans $\calD$ et les réseaux stables par Galois dans $V$. (Pour ce
dernier résultat, voir \cite{liu}.)

\paragraph{Les catégories $\ModphiN {S_1}$ et $\ModphiN {\tilde S_1}$}

On est maintenant tenté de réduire modulo $p$ les modules fortement
divisibles que l'on vient d'introduire pour obtenir une description des
représentations de la forme $T/pT$. On pose pour cela $S_1 = S/pS$. Les
structures supplémentaires sur $S$ passent au quotient pour définir des
structures analogues sur $S_1$.
On introduit la catégorie $\ModphiN {S_1}$ dont les objets sont la
donnée de :
\begin{itemize}
\item un module $\calM$ libre de type fini sur $S_1$ ;
\item un sous-module $\Fil^r \calM \subset \calM$ contenant $\Fil^r S_1
\: \calM$ ;
\item une application semi-linéaire $\phi_r : \Fil^r \calM \to \calM$
dont l'image engendre $\calM$ sur $S_1$ ;
\item une application $N : \calM \to \calM$ vérifiant la condition de
Leibniz, telle que $u^e N(\Fil^r \calM) \subset \Fil^r \calM$ et faisant
commuter le diagramme suivant :
$$\xymatrix @C=30pt {
\Fil^r \calM \ar[r]^-{\phi_r} \ar[d]_-{u^e N} & \calM \ar[d]^-{c N} \\
\Fil^r \calM \ar[r]^-{\phi_r} & \calM }$$
\end{itemize}
Si $er < p-1$, la catégorie $\ModphiN {S_1}$ est étudiée dans
\cite{caruso-crelle} : on montre en particulier qu'elle est abélienne et
artinienne et on décrit ses objets simples lorsque le corps résiduel $k$
est algébriquement clos. De plus, elle est également équipée d'un
foncteur exact $T_\st$ vers la catégorie des $\F_p$-représentations de
$G_K$.
Finalement, si $\calM$ est un module fortement divisible, le quotient
$\calM / p\calM$ est un objet de $\ModphiN {S_1}$ et si $T$ est le
réseau associé à $\calM$ par le foncteur $T_\st$, on a
$T_\st(\calM/p\calM) = T/pT$.

Pour les calculs de la section \ref{sec:exemple}, nous aurons besoin de
manipuler la catégorie $\ModphiN {\tilde S_1}$ : sa définition est
analogue à celle de $\ModphiN {S_1}$ hormis que $S_1$ est partout
remplacé par $\tilde S_1 = k[u]/u^{ep}$. Le morphisme d'anneaux $S_1 \to
\tilde S_1$ qui envoie $u$ sur $u$ et les puissances divisées
$\gamma_i(u^e)$ sur $0$ pour $i \geq p$ définit \emph{via} extension des
scalaires un foncteur $\ModphiN {S_1} \to \ModphiN {\tilde S_1}$. La
proposition 2.3.1 de \cite{caruso-crelle} assure que ce dernier est une
équivalence de catégories.

\section{Position relative des divers polygones}
\label{sec:position}

Nous prouvons à présent les théorèmes \ref{theo:hodge} et
\ref{theo:newton}. Nous commençons par quelques rappels sur les bases
adaptées qui jouent un rôle essentiel dans la démonstration du théorème
\ref{theo:hodge}.

\subsection{Notion de bases adaptées}

\subsubsection*{Modules libres sur les anneaux \og principaux \fg}

Nous isolons ici un lemme, conséquence facile du théorème de structure
des modules de type fini sur les anneaux principaux. Bien
qu'élémentaire, ce lemme peut être vu comme la clé de la démonstration
du résultat de comparaison entre polygone de Hodge et de l'inertie
modérée.

\begin{lemme}
\label{lem:divis}
Soit $A$ un anneau principal et $\mathfrak p$ un élément non nul
irréductible de $A$. On note encore $\mathfrak p \subset A$ l'idéal
principal maximal engendré par $\mathfrak p$.  Soient $r < N$ et $d$ des
entiers, et $M'$ un sous-$A$-module de $M = (A/\mathfrak p^N)^d$ tel que
$\mathfrak p^r A \subset M'$. Alors, il existe $(e_1, \ldots, e_d)$ une base
de $M$ (sur $A/\mathfrak p^N$) et des entiers $0 \leq n_1 \leq n_2 \leq
\cdots n_d \leq r$ tels que $M'$ soit le sous-module engendré par les
$\mathfrak p^{n_i}e_i$.

De plus, les $n_i$ sont uniquement déterminés (\emph{i.e.} ne dépendent
pas de la base $(e_1, \ldots, e_d)$) et peuvent s'obtenir de la façon
suivante. Soit $(x_1, \ldots, x_D)$ une famille génératrice de $M'$.
Soient $G$ la matrice $d \times D$ obtenue en écrivant en colonne les
composantes des vecteurs $x_i$ et $\hat G$ une matrice à coefficients
dans $A$ relevant $G$. Alors, pour tout $k \leq d$, le nombre $n_1 +
\cdots + n_k$ est la plus petite valuation $\mathfrak p$-adique d'un mineur
$k \times k$ de $\hat G$.
\end{lemme}

\begin{proof}
D'après la théorie des diviseurs élémentaires, il existe des matrices
inversibles $P$ et $Q$ (à coefficients dans $A$) telles que :
$$\hat G = P \cdot \pa{\raisebox{0.5\depth} {\xymatrix @C=3pt {
a_1 \ar@{.}[rrd] & & & 0 \ar@{.}[r] \ar@{.}[d] & 0 \ar@{.}[d] \\
& & a_d & 0 \ar@{.}[r] & 0 }}} \cdot Q$$
où $a_k$ est le $\pgcd$ des mineurs $k \times k$ de $\hat G$. De ceci,
il découle toutes les assertions du lemme à part l'unicité des $n_i$.
Pour cette dernière, on remarque que pour tout entier $n \leq r$ la
longueur du $A$-module $\mathfrak p^n M/(M' \cap \mathfrak p^n M)$ s'égalise
avec la somme des $\max(0, n_i-n)$, et que la connaissance de toutes ces
sommes permet de retrouver les $n_i$.
\end{proof}

\subsubsection*{Existence de bases adaptées}

Le lemme \ref{lem:divis} donne directement le résultat suivant, déjà
bien connu.

\begin{prop}
Soit $\calM$ un module filtré sur $S_{K_0}$ (resp. un objet
de $\ModphiN {S_1}$). Il existe $(e_1, \ldots, e_d)$ une base de $\calM$
et des entiers $n_1, \ldots, n_d$ compris entre $0$ et $r$ (resp. entre
$0$ et $er$) tels que :
$$\Fil^r \calM = \Fil^p S_{K_0} \: \calM + \sum_{i=0}^d E(u)^{n_i} e_i
S_{K_0}$$ $$\text{(resp. }\Fil^r \calM = \Fil^p S_1 \: \calM +
\sum_{i=0}^d u^{n_i} e_i S_1 \text{).}$$
De plus les entiers $n_i$ sont uniquement déterminés (à l'ordre près).
\end{prop}

\begin{proof}
On applique le lemme \ref{lem:divis} respectivement aux quotients
$S_{K_0}/ \Fil^p S_{K_0} \simeq K_0[u] / E(u)^p$ et $S_1/\Fil^p S_1
\simeq k[u]/u^{ep}$.
\end{proof}

\medskip

\noindent
Une base comme en fournit la proposition est appelé une \emph{base
adaptée} pour les entiers $n_1, \ldots, n_d$. On prendra garde au fait
qu'il n'est pas vrai que tout module fortement divisible admet une base
adaptée. Toutefois, c'est le cas si $r=1$ (et la démonstration repose
encore sur la théorie des diviseurs élémentaires, l'anneau principal qui
intervient étant ici $S / \Fil^1 S \simeq \O_K$).

\bigskip

L'égalité (\ref{eq:filt}) conduit directement à la proposition suivante :

\begin{prop}
\label{prop:hodge}
Soit $V$ une représentation semi-stable de $G_K$ et $\calD$ son
$S_{K_0}$-module filtré associé. Si $n_1, \ldots, n_d$ sont les entiers
qui apparaissent dans l'écriture d'une base adaptée de $\calD$, alors
les poids de Hodge-Tate de $V$ sont les $h_i = r - n_i$.
\end{prop}

\subsection{Preuve du théorème \ref{theo:hodge}}

On suppose $er < p-1$. On fixe $V$ une représentation semi-stable de
$G_K$ de dimension $d$. On note $\calD$ son $S_{K_0}$-module filtré
associé. D'après le résultat principal de \cite{breuil-invent}, il
existe $\calM \subset \calD$ un réseau fortement divisible. La
représentation galoisienne $T = T_\st (\calM)$ est un réseau dans $V$
stable par $G_K$ et on a $T/pT = T_\st(\calM / p\calM)$. Le but est de
comparer les poids de Hodge-Tate de $V$, notés $h_1 \leq \cdots \leq
h_d$, avec les poids de l'inertie modérée de $T/pT$, notés $i_1 \leq
\cdots \leq i_d$.

Notons $n_1 \geq \cdots \geq n_d$ (resp. $n'_1 \geq \cdots \geq n'_d$)
les entiers qui interviennent dans l'écriture d'une base adaptée de
$\calD$ (resp. de $\calM/p\calM$). Par la proposition \ref{prop:hodge},
on a la relation $h_k = r - n_k$ pour tout $k \in \{1, \ldots, d\}$.
Définissons les poids de Hodge-Tate de $\calM/p\calM$ comme les
rationnels $h'_k = r - \frac{n'_k} e$ et appelons polygone de Hodge de
$\calM/p\calM$ le polygone associé à ces nombres. La démonstration se
découpe alors naturellement en deux étapes : tout d'abord on montre que
le polygone de Hodge de $\calM/p\calM$ est au-dessus de celui de $V$
(avec même point terminal), puis qu'il est au-dessous du polygone de
l'inertie modérée (avec également même point terminal).

\subsubsection*{Première étape : comparaison entre polygones de Hodge}

On conserve les notations précédentes. Montrons tout d'abord que pour tout 
$k \in \{1, \ldots, d\}$, on a :
$$e (n_d + \ldots + n_k) \leq n'_d + \ldots + n'_k$$
et que l'égalité a lieu lorsque $k = 1$.

\medskip

Soit $(e_1, \ldots, e_d)$ une base adaptée de $\calM / p \calM$, donc
associée aux entiers $n'_1, \ldots, n'_d$. Notons $\hat x_i \in \Fil^r
\calM$ un relevé de $u^{n'_i} e_i$ et $x_i$ sa projection dans $\Fil^r
\calM / \Fil^p S \: \calM$. On vérifie facilement que la famille des
$x_i$ engendre $\Fil^r \calM / \Fil^p S\: \calM$.  Soit $G$ la matrice
carrée de taille $d$ obtenue en écrivant en colonne les coordonnées des
vecteurs $x_i$ dans une base quelconque (mais fixée) de $\calM / \Fil^p
S \: \calM$.  C'est une matrice à coefficients dans $W[u]/E(u)^p$.
Soient $\hat G$ une matrice à coefficients dans $W[u]$ relevant $G$ et
$\bar G$ la réduction de $\hat G$ dans l'anneau $k[u]$. 

Fixons un entier $k \in \{1, \ldots, d\}$. D'après le lemme
\ref{lem:divis}, $n = n_d + \ldots + n_k$ (resp. $n' = n'_d + \ldots +
n'_k$) est la plus petite valuation $E(u)$-adique (resp $u$-adique) d'un
mineur $(d+1-k) \times (d+1-k)$ de $\hat G$ (resp. de $\bar G$). Ainsi,
$E(u)^n$ divise certainement tous les mineurs de taile $d+1-k$ de $\hat
G$. Par réduction modulo $p$, il s'ensuit que $u^{en}$ divise tous les
mineurs de taille $d+1-k$ de $\bar G$. D'où $en \leq n'$ comme annoncé.

Il reste à montrer que l'égalité a lieu lorsque $k = 1$. Par hypothèse
$\Fil^r S \: \calM \subset \Fil^r \calM$. On en déduit qu'il existe une
matrice $\hat H$ à coefficients dans $W[u]$ telle que $\hat G \hat H
\equiv E(u)^r I \pmod{E(u)^p}$ (où $I$ désigne la matrice identité).
Comme $p-r \geq 1$, il existe une matrice $C$ à coefficients dans $W[u]$
telle que $\hat G \hat H = E(u)^r (I - E(u) C)$. En considérant les
déterminants, il vient :
$$\det (\hat G) \: \det (\hat H) = E(u)^{rd} \: \Delta$$
où $\Delta \in W[u]$ est congru à $1$ modulo $E(u)$. L'anneau $W[u]$
étant factoriel, $\det(\hat G)$ prend la forme $E(u)^n \delta$ où $n
\leq rd$ est un entier et où $\delta$ divise $\Delta$. En particulier
$E(u)$ ne divise pas $\delta$ et la valuation $E(u)$-adique de $\det
(\hat G)$ est $n$. Ainsi $n = n_1 + \cdots + n_d$. D'autre part, en
réduisant modulo $p$, on obtient $\det(\bar G) = u^{en} \bar \delta$ où
$\bar \delta$ divise un élément congru à $1$ modulo $u^e$. En
particulier sa valuation $u$-adique est nulle et donc celle de
$\det(\bar G)$ vaut $en$. Le résultat annoncé s'ensuit.

\bigskip

Il résulte de ceci le résultat de comparaison, objet de ce paragraphe.

\begin{lemme}
Avec les notations précédentes, le polygone de Hodge de $V$ est
au-dessous du polygone de Hodge de $\calM/p\calM$. De plus, ils ont même 
point d'arrivée.
\end{lemme}

\begin{proof}
En soustrayant $e(n_1 + \cdots + n_d) = n'_1 + \cdots + n'_d$ des deux
côtés de l'inégalité précédemment prouvée, on obtient $e(n_1 + \cdots +
n_k) \geq n'_1 + \cdots + n'_k$. Le résultat s'obtient alors en divisant
par $(-e)$, puis en ajoutant $kr$ à cette dernière inégalité. Bien
entendu, le cas d'égalité se traite par la même manipulation algébrique.
\end{proof}

\subsubsection*{Deuxième étape : polygone de Hodge de $\calM/p\calM$ et
polygone de l'inertie modérée}

Nous prouvons en fait de façon plus générale le lemme suivant :

\begin{lemme}
Soit $\calN$ un objet de $\ModphiN {S_1}$. Le polygone de Hodge de
$\calN$ est au-dessous du polygone de l'inertie modérée de
$T_\st(\calN)$. De plus, ils ont même point d'arrivée.
\end{lemme}

La catégorie $\ModphiN {S_1}$ étant abélienne et artinienne, il suffit
de prouver d'une part que le résultat est vrai pour les objets simples
de cette catégorie et d'autre part qu'il passe aux extensions. Comme les
poids de l'inertie modérée ne dépendent que de l'action du groupe
d'inertie et que les poids de Hodge-Tate sont invariants par extension
non ramifiée, on peut supposer que le corps résiduel $k$ est
algébriquement clos. Dans ce cas, les objets simples sont décrits dans
\cite{caruso-crelle} (théorème 4.3.2) et la représentation galoisienne
associée à ceux-ci est également calculée dans \emph{loc. cit.}
(théorème 5.2.2). On constate alors sans difficulté que les polygones
sont bien disposés comme on le souhaite. (En réalité, les polygones sont
mêmes confondus, ici.)

Traitons à présent le cas des extensions. Donnons-nous :
$$0 \to \calN' \to \calN \to \calN'' \to 0$$
une suite exacte dans $\ModphiN {S_1}$. Notons $a'_1 \leq \cdots \leq
a'_{d'}$ (resp. $b'_1 \leq \cdots \leq b'_{d'}$) les poids de Hodge-Tate
de $\calN'$ (resp. de l'inertie modérée de $T_\st(\calN')$) et $a''_1
\leq \cdots \leq a''_{d''}$ (resp. $b''_1 \leq \cdots \leq b''_{d''}$)
ceux correspondant à $\calN''$. \'Etant donné que le foncteur $T_\st$
est exact et que les poids de l'inertie modérée de $T_\st(\calN)$ ne
dépendent que des quotients de Jordan-Hölder de ladite représentation,
ils sont exactement les nombres $b'_1, \ldots, b'_{d'}, b''_1, \ldots,
b''_{d''}$.

\'Evaluons maintenant les poids de Hodge-Tate de $\calN$, que nous
notons $a_1 \leq a_2 \leq \cdots \leq a_d$ avec $d = d' + d''$. D'après
la preuve du lemme \ref{lem:divis}, pour tout entier $n \leq er$, ils
vérifient la relation :
$$\dim_k \frac{u^n \calN}{u^n \calN \cap \Fil^r \calN} = \sum_{i=1}^d
\max(0, e(r-a_i)-n)$$
et l'on dispose bien entendu de formules analogues pour $\calN'$ et
$\calN''$.  Considérons les deux applications :
$$\frac{u^n \calN'}{u^n \calN' \cap \Fil^r \calN'} \to \frac{u^n \calN}
{u^n \calN \cap \Fil^r \calN} \to \frac{u^n \calN''}{u^n \calN'' \cap 
\Fil^r \calN''}.$$
La seconde est clairement surjective et la première est injective : en
effet, si $x \in u^n \calN'$ s'envoie dans $\Fil^r \calN$, il est aussi
élément de $\Fil^r \calN'$ en vertu de la stricte compatibilité à la
filtration (voir corollaire 3.5.7 de \cite{caruso-crelle}). Il s'ensuit
que la dimension du terme central est supérieure à la somme des
dimensions des termes extrémaux. D'où on déduit l'inégalite suivante :

\begin{equation}
\label{eq:ineght}
\sum_{i=1}^d \max(0, e(r-a_i)-n) \geq \sum_{i=1}^{d'} \max(0, 
e(r-a'_i)-n) + \sum_{i=1}^{d''} \max(0, e(r-a''_i)-n).
\end{equation}
Il en résulte que pour tout $k \in \{1, \ldots, d\}$ la somme des $k$
plus petits entiers parmi les $a_i$ est inférieure à la somme des $k$
plus petits éléments parmi $a'_1, \ldots, a'_{d'}, a''_1, \ldots,
a''_{d''}$.

Par ailleurs, si $n = 0$, la suite :
$$0 \to \frac {\calN'} {\Fil^r \calN'} \to \frac \calN {\Fil^r \calN} \to
\frac {\calN''} {\Fil^r \calN''} \to 0$$
est exacte. En effet, il suffit d'après ce qui précède de vérifier
l'exactitude au milieu. Pour cela, on considère $x \in \calN$ qui
s'envoie sur un élément $y \in \Fil^r \calN''$ et il s'agit de montrer
que $x$ s'écrit comme la somme d'un élément de $\calN'$ et d'un élément
de $\Fil^r \calN$. Par la stricte compatiblité à la filtration,
l'application $\Fil^r \calN \to \Fil^r \calN''$ est surjective, d'où il
existe $x' \in \Fil^r \calN$ qui s'envoie également sur $y \in \Fil^r
\calN''$. Mais alors $x-x'$ est nul dans $\calN''$ et donc est élément
de $\calN$. L'exactitude en découle. De celle-ci, on déduit que lorsque
$n = 0$, l'inégalité (\ref{eq:ineght}) est en fait une égalité,
c'est-à-dire :
$$\sum_{i=1}^d a_i = \sum_{i=1}^{d'} a'_i + \sum_{i=1}^{d''} a''_i.$$
En regroupant tout ce qui précède, on s'aperçoit que l'on vient de
prouver que le polygone de Hodge de $\calN$ est situé au-dessous du
polygone associé aux nombres (retriés) $a'_1, \ldots a'_{d'}, a''_1,
\ldots, a''_{d''}$, et qu'ils ont même point terminaux. Pour conclure,
il ne reste donc plus qu'à prouver le lemme suivant.

\begin{lemme}
Supposons que le polygone associé aux nombres $a'_1, \ldots, a'_{d'}$
(resp.  $a''_1, \ldots, a''_{d'}$) soit au-dessous du polygone associé aux
nombres $b'_1, \ldots, b'_{d''}$ (resp. $b''_1, \ldots, b''_{d''}$) et
qu'ils aient même points d'arrivée.  Alors le polygone associe à $a'_1,
\ldots, a'_{d'}, a''_1, \ldots, a''_{d''}$ est au-dessous du polygone
associé à $b'_1, \ldots, b'_{d'}, b''_1, \ldots, b''_{d''}$, et ils ont
même point d'arrivée.
\end{lemme}

\begin{proof}
L'assertion sur les points d'arrivée est immédiate.

Pour le reste, on peut évidemment supposer que les $a'_i$, les $a''_i$,
les $b'_i$ et les $b''_i$ sont rangés par ordre croissant. Soit $k$
compris entre $1$ et $d$. La somme des $k$ plus petits nombres parmi
$a'_1, \ldots, a'_{d'}, a''_1, \ldots, a''_{d''}$ est égale à :
\begin{equation}
\label{eq:min}
\min_{1 \leq m \leq k} (a'_1 + \cdots + a'_m) + (a''_1 + \cdots + 
a''_{k-m})
\end{equation}
et, bien entendu, on dispose d'une formule analogue lorsque la lettre
$a$ est remplacée par $b$. La conclusion s'obtient en remarquant que
chacun des termes qui apparaît dans le minimum de (\ref{eq:min}) est
majoré par hypothèse par le terme correspondant où $a$ est remplacé par
$b$, et que cette majoration se transporte directement sur les minimas.
\end{proof}

\subsubsection*{Remarques}

\paragraph{Le cas $\mathbf{er \geq p-1}$}

Comme nous le disions dans l'introduction, le théorème \ref{theo:hodge}
peut être mis en défaut si on ôte l'hypothèse $er < p-1$. Voici un
contre-exemple très simple. Considérons $K_0 = \Q_p$ et $K =
K_0(\sqrt[p] 1)$ ; l'extension $K/K_0$ est totalement ramifiée de degré
$e = p-1$. Le caractère cyclotomique a pour seul poids de Hodge-Tate $1$
mais sa réduction modulo $p$ est triviale, de sorte que le poids de
l'inertie modérée est $0$. On remarque qu'ici on peut choisir $r = 1$,
de sorte que $er = p-1$ qui est la première situation dans laquelle
l'inégalité de l'énoncé est violée.

\paragraph{Que dire des représentations de Hodge-Tate ?}

Si $V$ n'est pas semi-stable, mais simplement de Hodge-Tate, l'énoncé du
théorème \ref{theo:hodge} a encore un sens, et on peut légitimement se
demander s'il est encore vrai dans ce contexte plus général. La réponse
est négative et, là encore, les contre-exemples sont aisés à produire.
On prend $K = K_0 = \Q_p$, et on fixe $\pi \in \bar K$ une racine 
$(p-1)$-ième de $p$. Définissons le caractère $G_K \to \Gal(K(\pi)/K) 
\to \Q_p^\star$, $\sigma \mapsto \frac{\sigma \pi} {\pi}$. Il correspond 
à une représentation
$p$-adique de dimension $1$ qui est de Hodge-Tate\footnote{Elle est même
\og potentiellement triviale \fg. En particulier, elle est
potentiellement cristalline.} et son seul poids de Hodge-Tate est nul.
Toutefois, modulo $p$, cette représentation s'identifie par construction
au caractère fondamental de Serre : son poids de l'inertie modérée est
donc $1$ et le résultat est mis en défaut.

\subsection{Preuve du théorème \ref{theo:newton}}
\label{subsec:newton}

On conserve notre entier positif $r$ vérifiant $er < p-1$. De plus comme
les invariants qui interviennent dans l'énoncé du théorème
\ref{theo:newton} ne changent pas après une extension non ramifiée, on
peut supposer que le corps résiduel $k$ est algébriquement clos.

Soit $V$ une représentation semi-stable dont on note les poids de
Hodge-Tate $h_1 \leq \cdots \leq h_d \leq r$ et les pentes de l'action
de Frobenius $n_1 \leq \cdots \leq n_d$. On suppose que les polygones de
Newton et de Hodge commencent par un segment de même pente, c'est-à-dire
que $h_1 = n_1 = s$. Notons $d' \geq 1$ le plus grand entier pour lequel
$n_{d'} = s$. \'Etant donné que le polygône de Newton est situé
au-dessus du polygône de Hodge, on a $n_k = h_k = s$ pour tout $k \leq
d'$.

Soit $D \in \MF^\fa_K(\varphi,N)$ le $(\varphi,N)$-module filtré de
Fontaine associé à $V$.  Appelons $D_s$ la partie de pente $s$ de $D$.
D'après la définition de $d'$, c'est un sous-espace de dimension $d'$
stable $\varphi$. Par ailleurs l'opérateur de monodromie doit envoyer
$D_s$ sur $D_{s-1}$ (partie de pente $s-1$) mais celle-ci est nulle
puisque $s$ est la plus petite pente. Ainsi $N = 0$ sur $D_s$, et $D_s$
est stable par $\varphi$ et $N$. Soit $D' \subset D_s$ un sous-espace de
dimension $n$ stable par $\varphi$ et $N$. Par les conditions de faible
admissibilité, on $t_H(D') \leq t_N(D') = ns$. Par ailleurs, $\Fil^t D_K
= D_K$ pour tout $t \leq s$ et donc $\Fil^t D'_K = D'_K$ pour tout $t \leq
s$. On en déduit qu'il n'y a pas de saut dans la filtration avant $s$ et
donc que $t_H(D') \geq ns$. Finalement $t_H(D') = t_N(D') = ns$, et il
s'ensuit que $D_s$ est un sous-objet de $D$ dans la catégorie
$\MF^\fa_K(\varphi,N)$ auquel il correspond une sous-représentation
$V_s$ de $V$.

Comme on a supposé le corps $k$ algébriquement clos, il existe une base
de $D_s$ qui diagonalise l'action du Frobenius avec des valeurs propres
toutes égales à $p^s$ (on rappelle que $s$ est entier puisque c'est un
poids de Hodge-Tate). Il en résulte que $V_s \simeq \Q_p(s)^{d'}$. Il
est alors immédiat de calculer les poids de l'inertie modérée sur $V_s$
: on trouve qu'ils sont égaux à $es$ comme attendu.

\medskip

Les considérations précédentes démontrent le théorème \ref{theo:newton}
pour la première pente. Pour la suite, on considère le quotient $D/D_s$
auquel on réapplique les arguments précédents. On continue comme cela
jusqu'à ce que les polygones se séparent, et cela clôt la preuve du
théorème.

\subsubsection*{Remarques}

Le résultat précédent n'est pas réellement satisfaisant. On aimerait par
exemple avoir une comparaison entre polygone de l'inertie modérée et
polygone de Newton qui ne repose sur aucune hypothèse. Bien entendu, la
question la plus simple que l'on puisse poser est la suivante :

\begin{question}
\label{quest:newton}
On suppose toujours que les poids de Hodge-Tate sont compris entre $0$
et $r$ (avec $er < p-1$). Est-il vrai que le polygone de l'inertie
modérée est situé entre le polygone de Hodge et celui de Newton ?
\end{question}

\noindent
On peut également s'interroger sur la positive relative du polygone
de Newton et celui de Hodge de la réduction modulo $p$ :

\begin{question}
\label{quest:newton2}
On suppose $r < p - 1$.
Soit $\calM$ un module fortement divisible. Est-il vrai que le polygone
de Hodge de $\calM/p\calM$ est situé dessous le polygone de Newton de
$D = \calM \otimes_S K_0$ ?
\end{question}

\noindent
L'avantage de la question \ref{quest:newton2} est qu'elle poser pour
tout $r < p-1$ : on entend par là qu'elle n'est semble-t-il pas
trivialement mise en défaut lorsque $er \geq p-1$ comme, nous l'avons
déjà expliqué, c'est le cas pour la question \ref{quest:newton}. En
revanche, elle propose un résultat un peu plus faible. Hélas, on dispose
de trop peu d'exemples pour se faire une idée claire de la réponse aux
questions précédentes. Lorsque $e = 1$, dans \cite{breuil-mezard},
Breuil et Mézard ont classifié complètement les représentations
semi-stables \emph{de dimension $2$} et ont calculé les poids de
l'inertie modérée pour chacune d'entre elles. Leurs résultats sont
résumés dans la partie 5 de \cite{breuil-azumino}, et montrent que la
réponse à la question \ref{quest:newton} (et donc aussi à la question
\ref{quest:newton2}) est affirmative dans ce cas. Toutefois, cela ne
doit pas être pris comme un indicateur très fort car dans la situation
étudiée ici, soit la représentation est cristalline et alors le polygone
de l'inertie modérée s'identifie avec le polygone de Hodge par les
résultats de Fontaine-Laffaille, soit la représentation n'est pas
cristalline, et alors le polygone de Newton est \og le plus haut
possible \fg.

\subsection*{Sans l'opérateur de monodromie}

La constatation initiale est la suivante : les objets polygone de Hodge,
polygone de Hodge de la réduction modulo $p$, polygone de Newton et
polygone de l'inertie modérée se calculent directement sur les modules
de la théorie de Breuil, indépendamment des représentations. Mieux
encore, leur calcul ne fait jamais intervenir l'opérateur de monodromie
$N$. Précisément, on a vu que le polygone de Hodge et le polygone de
Hodge de la réduction modulo $p$ s'obtiennent à partir des entiers qui
interviennent dans l'écriture d'une base adaptée des modules
correspondants. Le polygone de Newton d'un $S_{K_0}$-module filtré
$\calD$ est construit à partir des pentes de l'action du Frobenius sur
$\calD \otimes_{S_{K_0}} K_0$. Finalement, le polygone de l'inertie
modérée de $\calM$ peut s'obtenir en recollant les polygones de Hodge
des quotients de Jordan-Hölder\footnote{On peut montrer (voir
\cite{caruso-crelle}) qu'il existe une sous-catégorie commune à
$\ModphiN {S_1}$ à $\Modphi {S_1}$ (même définition que $\ModphiN {S_1}$
mais sans le $N$) qui contient tous les objects simples de ces deux
catégories. Ainsi, avec la définition que l'on donne, le polygone de
l'inertie modérée est le même qu'il soit calculé dans $\ModphiN {S_1}$
ou $\Modphi {S_1}$.} de $\calM/p\calM$, ce dernier n'était défini que
pour $er < p-1$. 

\medskip

Par ailleurs, on s'aperçoit que dans les énoncés des théorèmes
\ref{theo:hodge} et \ref{theo:newton}, il revenait au même de débuter
avec un module fortement divisible plutôt qu'une représentation
semi-stable. En effet, si l'on part d'une représentation semi-stable, on
peut trouver un module fortement divisible dans son $S_{K_0}$-module
filtré associé comme cela a déjà été expliqué. Et réciproquement, si
l'on part d'un module fortement divisible, le $S_{K_0}$-module filtré
obtenu en inversant $p$ correspond à un module filtré de Fontaine qui
est faiblement admissible (corollaire 2.1.4 de \cite{breuil-invent}), et
donc à une représentation galoisienne semi-stable.

\bigskip

Au vu de ces remarques, il est légitime de se demander si les théorèmes
\ref{theo:hodge} et \ref{theo:newton} s'étendent à une situation plus
générale où l'on débuterait avec un pseudo-module fortement divisible
(c'est-à-dire un module fortement divisible sans l'opérateur $N$). Pour
le théorème \ref{theo:hodge}, la réponse est affirmative et la
démonstration est textuellement la même que celle que nous avons déjà
donnée puisqu'elle ne fait aucunément intervenir l'opérateur de
monodromie. Les auteurs, par contre, ne savent pas ce qu'il en est pour
le théorème \ref{theo:newton}. Toutefois, il est possible de produire un
contre-exemple à la généralisation évidente de la question
\ref{quest:newton} à cette situation.

On suppose que $r = 2n$ est un entier pair, que $e = 1$ et que
l'uniformisante choisie est $p$ de sorte que $E(u) = u-p$. Soit le
pseudo-module fortement divisible défini par $\calM = S e_1 \oplus S
e_2$. On le munit de $\Fil^r \calM$ engendré par $(u-p)^n e_1$, $p e_1 +
(u-p)^n e_2$ et $\Fil^p S \: \calM$, et du Frobenius $\phi$ défini par
$\phi(e_1) = c^n p^n e_2$ et $\phi(e_2) = c^n p^n e_1 - p e_2$. Les
entiers qui apparaissent dans l'écriture d'une base adaptée de
$\calM/p\calM$ se calculent directement par le lemme \ref{lem:divis} et
valent $n$ et $n$. Le polygone de Hodge de $\calM/p\calM$ a donc pour
seule $n$, et comme le polygone de l'inertie modérée doit être situé
au-dessus, il a lui aussi pour seule pente $n$. Cependant, on calcule
facilement les pentes du Frobenius sur $\calM \otimes_S K_0$ : elles
valent $1$ et $r-1$. Le polygone de Newton est donc strictement
en-dessous celui de l'inertie modérée.

\section{Un (contre-)exemple dans le cas cristallin}
\label{sec:exemple}

Comme nous l'avons dit dans l'introduction, les travaux de Fontaine et
Laffaille assurent que dans le cas cristallin (\emph{i.e.} $N = 0$ sur
le $(\varphi,N)$-module filtré de Fontaine) non ramifié (\emph{i.e.} $e
= 1$), le polygone de l'inertie modérée s'identifie avec celui de Hodge.
Nous avons également déjà dit que Breuil et Mézard ont montré par un
calcul explicite que ce résultat simple est mis en défaut dès que la
représentation n'est plus supposée cristalline (même dans le cas $e=1$). 
Également par un calcul explicite, nous ne proposons de montrer dans
cette dernière section que l'identification entre polygone de l'inertie
modérée et polygone de Hodge n'est pas non plus satisfaite en général
pour les représentations cristallines dès que $e \geq 2$ (voir théorème
\ref{theo:red} et la remarque qui le suit). Le cas de Fontaine-Laffaille
apparaît donc, de ce point de vue, comme très isolé.

\medskip

Soient $n_1 \leq n_2$ des entiers positifs ou nuls tels que $e(n_1 +
n_2) < p-1$ (on prendra relativement rapidement $n_1 = n_2 = 1$, mais on
préfère commencer avec un peu plus de généralité). On pose $r = n_1 +
n_2$ et on définit pour tout $\calL \in K$, le $(\varphi,N)$-module
filtré suivant :
$$\left\{\begin{array}{l}
  D(\calL) = K_0 e_1 \oplus K_0 e_2 \\
  \varphi(e_1) = p^{n_1} e_1, \, \varphi(e_2) = p^{n_2} e_2 \\ 
  \Fil^1 D(\calL)_K = \cdots = \Fil^r D(\calL)_K = K(\calL e_1 + e_2) \\
  \Fil^0 D(\calL)_K = D(\calL)_K, \, \Fil^{r+1} D(\calL)_K = 0 \\
  N = 0 \end{array} \right.$$
Dans le cas où $n_1 = n_2$, ce module est admissible si, et seulement
si $\calL \not\in \Q_p$. Si $n_1 \neq n_2$, il est toujours admissible
sauf précisément dans le cas où $n_1 > 0$ et $\calL = 0$. On suppose
à partir de maintenant que $\calL$ est choisi de façon à ce que 
$D(\calL)$ soit admissible.

Le polygone de Hodge de $D(\calL)$ est directement visible sur la
description précédente : c'est celui qui a pour pentes $0$ et $r$. Notre
objectif désormais est de déterminer (au moins dans certains cas), le
polygone de l'inertie modérée associé à $D(\calL)$. Pour cela, nous
allons successivement déterminer le module filtré sur $S_{K_0}$ qui lui
correspond (sous-section \ref{subsec:SK0}), trouver un réseau fortement
divisible à l'intérieur de ce dernier (sous-section \ref{subsec:reseau})
et finalement réduire celui-ci modulo $p$ (sous-section
\ref{subsec:redmodp}). 

\subsection{Calcul du module filtré sur $S_{K_0}$}
\label{subsec:SK0}

Par définition le $S_{K_0}$-module filtré associé à $D(\calL)$ est
$\calD(\calL) = S_{K_0} \otimes_{K_0} D(\calL)$ muni du Frobenius et
de l'opérateur de monodromie obtenus par extension des scalaires. On
a donc directement
$$\left\{\begin{array}{l}
\calD(\calL) = S_{K_0} e_1 \oplus S_{K_0} e_2 \\
\phi(e_1) = p^{n_1} e_1, \, \phi(e_2) = p^{n_2} e_2 \\ 
N(e_1) = 0, \, N(e_2) = 0 \\
\end{array} \right.$$
où dans un léger abus de notation, on continue à noter $e_1$ et $e_2$
les éléments $1 \otimes e_1$ et $1 \otimes e_2$. Déterminer la forme
de la filtration demande par contre un peu de travail. Considérons,
l'application :
$$\begin{array}{rcl}
T_\pi : K_0[u]/E(u)^r & \longrightarrow & K^r \\
P & \mapsto & (P(\pi), P'(\pi), \ldots, P^{(r-1)}(\pi))
\end{array}$$
où $P^{(i)}$ désigne la dérivée $i$-ième du polynôme $P$. Il est clair
que $T_\pi$ est $K_0$-linéaire et on vérifie facilement qu'elle est
injective. Comme les espaces de départ et d'arrivée sont tous les deux
de dimension $er$ sur $K_0$, $T_\pi$ est une bijection. Notons $\calL_r$
l'unique polynôme de degré $< er$ dont l'image par $T_\pi$ est $(\calL,
0, \ldots, 0)$. Bien entendu $\calL_r$ peut également être vu comme un
élément de $S_{K_0}$.

\begin{prop}
On a $\Fil^r \calD(\calL) = (\calL_r e_1 + e_2) S_{K_0} + \Fil^r
S_{K_0}$.
\end{prop}

\begin{proof}
Après avoir vérifié quelques compatibilités, la proposition peut
s'obtenir comme une application de la formule donnée dans la remarque
finale de \cite{breuil-griff} (celle qui suit la preuve de la
proposition A.4).
Nous donnons malgré tout ci-dessous une démonstration plus directe qui
n'utilise que la définition de la filtration. Précisément, nous allons
prouver par récurrence sur $s$ que $\Fil^s \calD(\calL) = (\calL_r e_1 +
e_2) S_{K_0} + \Fil^s S_{K_0}$ pour tout $s$ compris entre $0$ et $r$.
La proposition s'ensuivra directement.

L'initialisation de la récurrence est immédiate. Fixons à présent un
entier $s < r$ et supposons que l'on ait réussi à montrer $\Fil^s
\calD(\calL) = (\calL_{r+1} e_1 + e_2) S_{K_0} + \Fil^s S_{K_0}$.  Par
définition
$$\Fil^{s+1} \calD(\calL) = \{ x \in \calD \, / \, N(x) \in \Fil^s
\calD(\calL), \, f_\pi(x) \in \Fil^t D(\calL)\}.$$
La dernière condition \og $f_\pi(x) \in \Fil^t D(\calL)$ \fg\ étant 
équivalente à $x \in \Fil^1 \calD(\calL)$, on obtient en utilisant la 
décroissante de la filtration et $N(\Fil^{s+1} \calD(\calL)) \subset 
\Fil^s \calD(\calL)$, la suite d'implications suivante :
$$\begin{array}{rcccl}
(x \in \Fil^{s+1} \calD(\calL)) &\Rightarrow &
(x \in \Fil^s \calD(\calL), \, N(x) \in \Fil^s \calD(\calL)) \\
&\Rightarrow& (x \in \Fil^1 \calD(\calL), \, N(x) \in \Fil^s \calD(\calL)) 
& \Rightarrow & (x \in \Fil^{s+1} \calD(\calL)).
\end{array}$$
Toutes ces implications sont donc des équivalences, d'où on déduit que
$\Fil^{s+1} \calD(\calL)$ s'identifie à l'ensemble des $x \in \Fil^s 
\calD(\calL)$ pour lesquels $N(x) \in \Fil^s \calD(\calL)$.
Considérons maintenant $x \in \Fil^s \calD(\calL)$. D'après l'hypothèse 
de récurrence, il existe des éléments $A$ et $B$ dans $S_{K_0}$ tels que
$x \equiv (A \calL_t + B E(u)^s) \hat e_1 + A \hat e_2 \pmod {\Fil^{s+1}
S_{K_0} \calD(\calL)}$. Un calcul direct donne :
\begin{eqnarray*}
N(x) & \equiv & (N(A) \calL_r + A N(\calL_r) + t B E(u)^{s-1}
N(E(u))) \hat e_1 + N(A) \hat e_2 \\ 
& \equiv & (N(A) \calL_r + t B E(u)^{s-1} N(E(u))) \hat e_1 + N(A) \hat
e_2 \pmod {\Fil^s S_{K_0} \calD(\calL)}
\end{eqnarray*}
la dernière congruence s'obtenant après avoir remarqué que $N(\calL_r) =
-u \calL_r'(u)$ est multiple de $E(u)^s$ étant donné que par définition
toutes ses dérivées d'ordre $< r-1$ (et donc \emph{a fortiori} d'ordre
$< s$) s'annulent. Ainsi, en utilisant à nouveau l'hypothèse de
récurrence, on trouve que $x \in \Fil^{s+1} \calD(\calL)$ si, et
seulement si $t B E(u)^{s-1} N(E(u)) \in \Fil^s S_{K_0}$, \emph{i.e.} $B
N(E(u)) \in \Fil^1 S_{K_0}$. Comme $S_{K_0} /\Fil^1 S_{K_0} \simeq K$
est un corps et que $N(E(u))$ ne s'annule pas dans ce quotient, cela
équivant encore à $B \in \Fil^1 S_{K_0}$, à partir de quoi l'hérédité
s'obtient directement.
\end{proof}

\subsection{Calcul d'un réseau fortement divisible}
\label{subsec:reseau}

Nous cherchons maintenant à construire un réseau fortement divisible à
l'intérieur de $\calD(\calL)$, essayant par là-même de généraliser
l'exemple 2.2.2.(4) de \cite{breuil-azumino} (qui est le point de départ
de tout notre calcul). En réalité, nous n'allons y parvenir
(probablement par manque de courage) que pour $n_1 = n_2 = 1$, ce qui
sera déjà suffisant pour mettre en défaut l'analogue du résultat de
Fontaine-Laffaille expliqué dans l'introduction de cette section. Nous
commençons toutefois par donner un outil pour construire des réseaux
fortement divisibles qui s'applique à tous $n_1$ et $n_2$, et donc ---
nous l'espérons --- pourra être utilisé fructueusement dans de futures
références.

\subsubsection{Un outil pour obtenir des réseaux fortement divisibles}

\begin{deftn}
Pour tout élément $X \in S_{K_0}$, on appelle \emph{$r$-ième troncation}
de $X$ et on note $\tr_r(X)$ l'unique polynôme en $u$ à coefficients
dans $K_0$ de degré $< er$ congru à $X$ modulo $\Fil^r S_{K_0}$.
\end{deftn}

Le lemme suivant réunit quelques propriétés immédiates de l'application
de troncation.

\begin{lemme}
\label{lem:tronc}
\begin{enumerate} 
\item L'application $\tr_r$ est $K_0$-linéaire.
\item Pour tous $X$ et $Y$ dans $S_{K_0}$, on a l'équivalence :
$$X \equiv Y \pmod {\Fil^r S_{K_0}} \quad \Longleftrightarrow \quad
\tr_r(X) = \tr_r(Y).$$
\item Pour tout $X \in S$, on a $\phi(X) \equiv \phi \circ
\tr_r (X) \pmod {p^r S}$.
\end{enumerate}
\end{lemme}

\begin{proof}
Les deux premières propriétés sont évidentes. Pour la dernière, il
suffit de remarquer que $X - \tr_r(X) \in S \cap \Fil^r S_{K_0} =
\Fil^r S$ et que $\phi(\Fil^r S) \subset p^r S$.
\end{proof}

\begin{prop}
\label{prop:condsuf}
On suppose donnés un entier relatif $n$ et un élément $Z \in p^{n-n_1} 
S$ satisfaisant l'implication suivante :
$$\H(\calL) : \,
\left\{\!\!\begin{array}{l}
A \in S \\
A (Z - \calL_r p^{n_2 - n_1}) \in p^n S + \Fil^r S_{K_0}
\end{array} \right.
\,\, \Longrightarrow \,\,\,
\left\{\!\!\begin{array}{l}
\phi(A) \in p^{n_1} S \\
\phi \circ \tr_r (A \calL_r) \equiv \phi(A) Z \pmod {p^{n+n_1} S}
\end{array} \right.$$
Alors le $S$-module engendré par $Z e_1 + p^{n_2-n_1} e_2$ et $p^n e_1$ 
est un réseau fortement divisible dans $\calD(\calL)$.
\end{prop}

\begin{proof}
Définissons $f_1 = Z e_1 + p^{n_2-n_1} e_2$ et $f_2 = p^n e_1$. Le
déterminant de la matrice $\Big(\,\begin{matrix} p^n & Z \\ 0 &
p^{n_2 -n_1} \end{matrix}\,\Big)$ étant manifestement une puissance de
$p$, il est clair que le $S$-module $\calM = S f_1 + S f_2$ est libre et
qu'il vérifie $\calM \otimes_S S_{K_0} = \calD$. Il ne reste donc qu'à
montrer que $\phi(\Fil^r \calM) \subset p^r \calM$ où, par définition,
$\Fil^r \calM = \calM \cap \Fil^r \calD$.
Soit $x \in \Fil^r \calM$. Alors $x \in \calM$, d'où on déduit qu'il
existe $A$ et $B$ dans $S$ tels que $x = A f_1 + B f_2$. Cette dernière 
égalité se réécrit encore :
$$x = A p^{n_2 - n_1} (\calL_r e_1 + e_2) + (p^n B + A Z - A \calL_r
p^{n_2-n_1}) e_1$$
ce qui, combiné à l'hypothèse \og $x \in \Fil^r \calD$ \fg\ et à la
définition de $\Fil^r \calD$, montre que
\begin{equation}
\label{eq:condfilr}
p^n B + A Z - A \calL_r p^{n_2-n_1} \in \Fil^r S_{K_0}.
\end{equation}
En particulier $A (Z - \calL_r p^{n_2-n_1}) \in p^n S + \Fil^r S_{K_0}$ et
le prémisse de $\H(\calL)$ est satisfait. Ainsi, par hypothèse,
$p^{n_1}$ divise $\phi(A)$ et $\phi \circ \tr_r (A \calL_r) \equiv
\phi(A) Z \pmod {p^{n + n_1} S}$. Rappelons que nous souhai\-tons obtenir
$\phi(x) \in p^r \calM$. Les égalités
\begin{eqnarray*}
\phi(x) & = & p^{n_1} \phi(AZ) e_1 + p^{2n_2 - n_1} \phi(A) e_2 + 
p^{n+n_1} \phi(B) e_1 \\
& = & p^{n_2} \phi(A) f_1 + \big[ p^{n_1-n} \phi(AZ) + p^{n_1} \phi(B) 
- p^{n_2-n} \phi(A) Z \big] f_2
\end{eqnarray*}
assurent que cela équivaut à démontrer que $p^r$ divise à la fois
$p^{n_2} \phi(A)$ et $p^{n_1-n} \phi(AZ) + p^{n_1} \phi(B) - p^{n_2-n}
\phi(A) Z$. Pour le premier, c'est immédiat puisque l'on sait que
$p^{n_1}$ divise $\phi(A)$. Pour le second, on remarque dans un premier
temps que la condition (\ref{eq:condfilr}) entraîne grâce aux sorites 
du lemme \ref{lem:tronc} :
$$p^n \: \tr_r(B) + \tr_r(AZ) = p^{n_2-n_1} \tr_r(A \calL_r)$$
puis :
$$p^n \phi(B) + \phi(AZ) \equiv p^{n_2-n_1} \phi \circ \tr_r(A \calL_r)
\pmod {p^{r+n-n_1} S}$$
(on rappelle que $Z$ est dans $p^{n-n_1} S$ par hypothèse).
On est donc ramené à prouver que $p^{n_2-n} \phi \circ \tr_r(A \calL_r)
\equiv p^{n_2-n} \phi(A) Z \pmod {p^r S}$, ce qui s'obtient directement
en multipliant par $p^{n_2-n}$ la congruence $\phi \circ \tr_r(A
\calL_r) \equiv \phi(A) Z \pmod {p^{n_1+n} S}$.
\end{proof}

\subsubsection{Le cas $n_1 = n_2 = 1$}

\emph{À partir de maintenant on suppose $n_1 = n_2 = 1$ (et donc $r=2$
et $p > 2e+1$).} Le but de ce paragraphe est de construire des réseaux
fortement divisibles à l'intérieur des $\calD(\calL)$. On commence par
un lemme qui va nous permettre de réduire notre étude à certains $\calL$
particuliers pour lesquels la proposition \ref{prop:condsuf} pourra être
appliquée.

\begin{lemme}
\label{lem:isoms}
Pour tous $\calL \in K$ et $a \in \Q_p$, les $S_{K_0}$-modules filtrés
$\calD(\calL)$ et $\calD(\calL+a)$ sont isomorphes.

Pour tous $\calL \in K$ et $n \in \Z$, les $S_{K_0}$-modules filtrés
$\calD(\calL)$ et $\calD(p^n \calL)$ sont isomorphes.
\end{lemme}

\begin{proof}
Pour la première assertion, l'isomorphisme est donné par $e_1 \mapsto
e_1$, $e_2 \mapsto a e_1 + e_2$. Pour la seconde assertion, il est
donné par $e_1 \mapsto p^n e_1$, $e_2 \mapsto e_2$.
\end{proof}

Ainsi, pour le calcul que l'on souhaite faire, on peut sans perte de
généralité remplacer $\calL$ par $p^n (\calL + a)$ avec $a \in \Q_p$ et
$n \in \Z$. Se faisant, il est facile de voir que l'on peut se ramener
à l'un de deux cas suivants (on rappelle que $\calL$ est supposé ne pas
appartenir à $\Q_p$) :
\begin{itemize}
\item[(i)] $v_p(\calL) = 0$ et l'image de $\calL$ dans le corps résiduel
n'appartient pas au sous-corps premier
\item[(ii)] $0 < v_p(\calL) < 1$.
\end{itemize}

\bigskip

Soit $L_0 \in K_0[u]$ l'unique polynôme de degré $< e$ tel que $L_0(\pi)
= \calL$. Comme l'on a supposé $0 \leq v_p(\calL) < 1$, $L_0$ est à
coefficients dans $W$. Par ailleurs son terme constant $\lambda$ est
inversible dans le cas (i) et multiple de $p$ dans le cas (ii). En
outre, dans le cas (i), l'hypothèse indique que $\mu = \phi(\lambda) -
\lambda$ est aussi inversible dans $W$. Dans le cas (ii), évidemment
$\mu$ est multiple de $p$. En utilisant la définition
de $\calL_2$, on montre qu'il s'écrit :
$$\calL_2 = L_0 + \frac 1 p L_1 E(u)$$ 
où $L_1$ est l'unique polynôme de degré $<e$ à coefficients dans $K_0$
tel que $L_1(\pi) = \frac{p L_0'(\pi)} {E'(\pi)}$. Comme $E(u)$ est un
polynôme d'Eisenstein de degré $e<p$, la valuation de $E'(\pi)$ est $1 -
\frac 1 e$. Ainsi $L_1(\pi)$ est divisible par $\pi$, d'où on déduit
$L_1 \in u S + p S$.

\begin{lemme}
\label{lem:petitt}
Il existe un élément $t \in S$ tel que :
\begin{itemize}
\item[$\bullet$] $(\phi(\lambda) - L_0) t \equiv L_1 \pmod {\Fil^1 S}$ ;
\item[$\bullet$] $1 + c \phi(t)$ est inversible dans $S$ (on rappelle 
que $c = \frac{\phi(E(u))} p$) ;
\item[$\bullet$] dans le cas (i), $t \in u S + pS$.
\end{itemize}
\end{lemme}

\begin{proof}
On traite séparément les deux cas. Dans le cas (i), on note que le terme
constant de $\phi(\lambda) - L_0$ est $\phi(\lambda) - \lambda = \mu$,
inversible dans $W$. Ainsi $\phi(\lambda) - L_0$ est lui-même inversible
dans $S$, et on pose $t = \frac {L_1}{\phi(\lambda) - L_0}$. Il vérifie
à l'évidence la première et la troisième condition. La seconde, quant à
elle, résulte directement de la troisième.

Passons maintenant au cas (ii). On choisit pout $t$ le polynôme à
coefficients dans $K_0$ de degré $< e$ tel que $t(\pi) =
\frac{L_1(\pi)}{\phi(\lambda) - L_0(\pi)}$. Il s'agit donc de montrer
d'une part que $t$ a ses coefficients dans $W$, et d'autre par la
deuxième condition du lemme (la première étant immédiate par
construction). Posons pour cela $j = e v_p(\calL)$. C'est un entier
strictement compris entre $0$ et $e$ et $L_0$ s'écrit :
$$L_0(u) = \ell_0 + \ell_1 u + \cdots + \ell_{e-1} u^{e-1}$$ 
où les $\ell_i$ sont des éléments de $W$ tels que $p$ divise $\ell_0,
\ldots, \ell_{j-1}$ tandis que $\ell_j$ est inversible. À partir de
cela, on obtient facilement les congruences $L_0(\pi) \equiv \ell_j
\pi^j \pmod{\pi^{j+1}}$ et $L'_0(\pi) \equiv j \ell_j \pi^{j-1}
\pmod{\pi^j}$. En faisant des manipulations analogues avec le polynôme
$E(u)$, il en ressort que $E'(\pi) \equiv e \pi^{e-1} \pmod {\pi^e}$ et
$p \equiv -\frac {\pi^e}{c_0} \pmod {\pi^{e+1}}$ si $p c_0$ est le
coefficient constant de $E(u)$. Ainsi trouve-t-on :
$$t(\pi) = \frac{L_1(\pi)}{\phi(\lambda) - L_0(\pi)} = \frac {-p 
L'_0(\pi)} {E'(\pi) (\phi(\lambda) - L_0(\pi))} \equiv - \frac j {e c_0} 
\pmod \pi$$
d'où il résulte que $t$ est à coefficients dans $W$ et que son 
coefficient constant est congru à $- \frac j {e c_0}$ modulo $p$. Par
suite, le coefficient constant de $1 + c \phi(t)$ est congru à $1 - 
\frac j e$ modulo $p$ et donc, en particulier, inversible dans $W$ 
(car $0 < j < e$). Il s'ensuit que $1 + c \phi(t)$ est inversible
dans $S$ comme annoncé.
\end{proof}

On fixe maintenant un élément $t \in S$ satisfaisant les conditions du
lemme précédent et on pose 
\begin{equation}
\label{eq:defZ}
Z = \frac{\phi(L_0) + c \phi(t \phi(\lambda))} { 1 + c \phi(t) } \in S.
\end{equation}

\begin{lemme}
\label{lem:Z}
On a $Z - \phi(\lambda) \in p S + \Fil^2 S$.
\end{lemme}

\begin{proof}
Un calcul donne :
$$Z - \phi(\lambda) = \frac{\phi(L_0 - \lambda) + c \phi(t \mu)}
{1 + c \phi(t)}.$$
Il suffit donc de montrer que $\phi(L_0 - \lambda) + c \phi(t \mu) \in p
S + \Fil^2 S$. Comme, par définition, $\lambda$ est le terme constant de
$L_0$, on a $L_0 - \lambda$ multiple de $u$ et donc $\phi(L_0 -
\lambda)$ est divisible par $u^p$. Or $u^p = u^{p-2e} u^{2e} \equiv
u^{p-2e} E(u)^2 \pmod {pS}$, d'où $u^p \in p S + \Fil^2 S$. Il en
résulte que $\phi(L_0 - \lambda)$ est toujours dans $p S + \Fil^2 S$.
Ainsi, pour terminer la preuve, il suffit de justifier que $\phi(t\mu)$
est lui aussi dans $p S + \Fil^2 S$. C'est clair si on est dans le cas
(ii) puisqu'alors $\mu$ est multiple de $p$. Si, au contraire, on est
dans le cas (i), le lemme \ref{lem:petitt} nous dit que $t \in uS + pS$
et on conclut comme précédemment.
\end{proof}

Notre but désormais est de montrer que le couple $(1,Z)$ satisfait
l'hypothèse $\H(\calL)$. Si on y parvient, par application de la
proposition \ref{prop:condsuf}, on aura bien réussi à construire un
réseau fortement divisible dans $\calD(\calL)$. On considère un élement
$A \in S$ satisfaisant le prémisse de $\H(\calL)$ (et on souhaite bien
évidemment montrer la conclusion). On observe qu'ajouter à $A$ un
élément de $\Fil^2 S$ ne modifie pas la véracité de $\H(\calL)$. Ainsi,
on peut supposer que $A = \tr_2(A)$ et par suite qu'il s'écrit sous la
forme $A = p A_0 + A_1 E(u)$ où $A_0$ (resp. $A_1$) est un polynôme à
coefficients dans $\frac 1 p W$ (resp. $W$) de degré $< e$.

\begin{lemme}
\label{lem:A0}
Avec les notations précédentes, $A_0$ est à coefficients dans $W$.
\end{lemme}

\begin{proof}
Encore une fois, on traite séparément les cas (i) et (ii).

Dans le cas (i), on remarque que $A(\calL_2 - Z)$ appartient
simultanément à $pS + \Fil^2 S_{K_0}$ et $\frac 1 p S$ (puisque $\calL_2
\in \frac 1 p S$). Ainsi $A(\calL_2 - Z) \in pS + \frac 1 p \Fil^2 S$ et
en appliquant $\phi$ on obtient $\phi(A) \phi(\calL_2 - Z) \in pS$. Or,
en utilisant le lemme \ref{lem:Z}, on vérifie facilement que
$\phi(\calL_2 - Z) = \phi(L_0) + c \phi(L_1) - \phi(Z)$ a un terme
constant congru à $-\phi(\mu)$ modulo $p$, et donc qu'il est inversible
dans $W$. Ainsi $\phi(\calL_2 - Z)$ est inversible dans $S$ et, de
$\phi(A) \phi(\calL_2 -Z) \in pS$, on déduit $\phi(A) \in pS$. Comme
$\phi(A) = p \phi(A_0) + p c \phi(A_1)$, il suit $\phi(A_0) \in S$ à
partir de quoi le lemme s'obtient facilement.

Dans le cas (ii), on a $A Z \in p S + \Fil^2 S_{K_0}$ et donc le
prémisse de $\H(\calL)$ s'écrit ici $A \calL_2 \in p S + \Fil^2
S_{K_0}$. Soit $X$ le polynôme à coefficients dans $W$ de degré $<e$ tel
que $X(\pi) = \frac p \calL \in W$. Alors, $\frac {\calL_2 X} p - 1$
s'annule en $\pi$, ce qui permet d'écrire 
\begin{equation}
\label{eq:defY}
\frac{\calL_2 X} p - 1 \equiv \frac Y p E(u) \pmod{\Fil^2 S_{K_0}}
\end{equation}
pour un certain polynôme $Y$ à coefficients dans $K_0$, uniquement
déterminé si on impose en outre $\deg Y < e$ (ce que nous ferons par la
suite). En dérivant (\ref{eq:defY}) et en évaluant en $\pi$, on obtient
$Y(\pi) = \frac{p X'(\pi)} {X(\pi) E'(\pi)}$. Un argument analogue à
celui conduit lors de la preuve du lemme \ref{lem:petitt} permet
d'accéder à la valuation de $Y(\pi)$ : on trouve $v_p(Y(\pi)) = 0$.
Ainsi $Y$ est à coefficients dans $W$ et son terme constant est
inversible dans cet anneau. On en déduit que $Y$ est inversible dans
$S$. Par ailleurs, on a
$$A_0 Y E(u) \equiv A \left( \frac{\calL_2 X} p - 1 \right) = \frac {A
\calL_2} p \: X - A
\pmod {\Fil^2 S_{K_0}}$$
d'où on trouve $A_0 Y E(u) \in S + \Fil^2 S_{K_0}$ puis $A_0 Y \in S +
\Fil^1 S_{K_0}$. De l'inversibilité de $Y$, on déduit enfin $A_0 \in S + 
\Fil^1 S_{K_0}$ à partir de quoi la conclusion est immédiate.
\end{proof}

En combinant les lemmes \ref{lem:tronc} et \ref{lem:A0}, il vient 
\begin{equation}
\label{eq:phiA}
\phi(A) \equiv \phi \circ \tr_2(A) = p \phi(A_0) + p c \phi(A_1) 
\pmod {p^2 S}.
\end{equation}
En particulier $\phi(A) \in pS$ (ce qui avait déjà été vu dans la preuve
du lemme \ref{lem:A0} dans le cas (i)). Il ne reste donc plus qu'à
démontrer que $p^2$ divise 
\begin{equation}
\label{eq:Delta}
\Delta = \phi \circ \tr_2(A \calL_2) - \phi(A) Z.
\end{equation}
En appliquant $\phi \circ \tr_2$ à la congruence
\begin{equation}
\label{eq:AL2}
A \calL_2 \equiv p A_0 L_0 + (A_0 L_1 + A_1 L_0) E(u) \pmod {\Fil^2
S_{K_0}}
\end{equation}
et en utilisant à nouveau le lemme \ref{lem:tronc}, on déduit :
\begin{equation}
\label{eq:phitr2}
\phi \circ \tr_2(A \calL_2) \equiv p \phi(A_0 L_0) + p c\phi(A_0 L_1 +
A_1 L_0) \pmod {p^2 S}.
\end{equation}
En remplaçant dans (\ref{eq:Delta}) $\phi \circ \tr_2(A \calL_2)$,
$\phi(A)$ et $Z$ par les expressions données respectivement par les
équations (\ref{eq:phitr2}), (\ref{eq:phiA}) et (\ref{eq:defZ}), on
obtient après un calcul un peu laborieux :
$$\frac{\Delta} p \equiv \frac{c \phi(A_0)}{1+c \phi(t)} \: \phi\big(
L_1 + L_0 t - \phi(\lambda) t\big) + \frac{c^2 \phi(t)}{1+c \phi(t)}
\: \phi\big(A_0 L_1 + A_1 L_0 - \phi(\lambda) A_1\big) \pmod {pS}.$$
Le premier terme du membre de droite est divisible par $p$ étant donné
que, par la première condition du lemme \ref{lem:petitt}, $L_1 + L_0 t - 
\phi(\lambda) t \in \Fil^1 S$. Il ne reste donc plus qu'à prouver qu'il
en est de même du second terme. Pour cela, on remarque que le lemme 
\ref{lem:Z} et la congruence (\ref{eq:AL2}) entraînent ensemble
\begin{equation} \label{eq:calA}
A(\calL_2 - Z) \equiv E(u) \big( A_0 L_1 + A_1 L_0 - A_1 
\phi(\lambda) \big) \pmod {p S + \Fil^2 S_{K_0}}.
\end{equation}
Ainsi, en utilisant l'hypothèse du prémisse de $\H(\calL)$, il vient
$E(u) ( A_0 L_1 + A_1 L_0 - A_1 \phi(\lambda) ) \in p S + \Fil^2 S$ puis
$A_0 L_1 + A_1 L_0 - A_1 \phi(\lambda) \in p S + \Fil^1 S$. En
appliquant $\phi$, on obtient finalement $\phi(A_0 L_1 + A_1 L_0 - A_1
\phi(\lambda) ) \in pS$ comme souhaité.

\subsection{Réduction modulo $p$ et poids de l'inertie modérée}
\label{subsec:redmodp}

On conserve les hypothèses et les notations de la partie précédente :
notamment, on continue de supposer que $\calL$ relève soit du cas (i),
soit du cas (ii) (on rappelle que, grâce au lemme \ref{lem:isoms}, on
peut toujours se ramener à l'un de ces deux cas), et en particulier donc
que $0 \leq v_p(\calL) < 1$. Soit $\calM = \calM(\calL)$ le réseau
fortement divisible de $\calD(\calL)$ qui est donné par la proposition
\ref{prop:condsuf}, et soit $T = T_{\st}(\calM)$ le $\Z_p$-réseau de
$V(\calL)$ correspondant. Posons $\barM = \calM/p\calM$, et pour tout $m
\in \calM$, notons $\overline{m}$ son image dans $\barM$. De même, si
$s$ est un élément de $S$, définissons $\overline{s}$ comme l'image de
$s$ dans $\tilde S_1 = k[u]/u^{ep}$ (\emph{cf} le dernier alinéa du
paragraphe \ref{subsec:torsion}). L'objectif de cette sous-section est
de démontrer le théorème suivant.

\begin{theo} 
\label{theo:red}
Le polygone de Hodge de $\barM$ a pour pentes $v_p(\calL)$ et $2 -
v_p(\calL)$.
De plus, si $v = v_p(t(\pi))$ (où on rappelle que $t \in S$ est un
élément satisfaisant les condition du lemme \ref{lem:petitt}), on a
\begin{enumerate}
\item si $0 \leq v < 1$, alors $T/pT$ est réductible et les pentes de
son polygone de l'inertie modérée sont $1-v$ et $1+v$ ;
\item si $v \geq 1$, alors $T/pT$ est irréductible et les pentes de
son polygone de l'inertie modérée sont $0$ et $2$ (\emph{i.e.} le
polygone de l'inertie modérée est confondu avec le polygône de Hodge).
\end{enumerate}
\end{theo}

\noindent
{\it Remarque.}
Lorsque $e=1$, on a $t(\pi) = 0$, et donc $T/pT$ est toujours
irréductible et son polygone de l'inertie modérée a pour pentes $0$ et
$2$. Si, au contraire, $e > 1$, alors tous les couples $(i_1', i_2')$
d'éléments de $\frac 1 e \N$ vérifiant $i_1' \leq i_2'$ et $i_1' + i_2'
= 2$ peuvent apparaître comme pentes du polygone de l'inertie modérée.
En effet, si $x$ est n'importe quel élément de $\O_{K_0}^{\times}$ dont
la réduction modulo $p$ n'appartient pas au sous-corps premier, et si $0
< j < e$ est un entier, alors les paramètres $\calL = \pi^j$, $x +
\pi^j$ et $x$ conduisent respectivement à $v = 0$, $\frac{j}{e}$ et
$\infty$.

\bigskip

Avant de commencer la démonstration, introduisons quelques notations
supplémentaires. Soient $U$ et $V$ les uniques polynômes à coefficients
dans $W$ de degré $< e$ tels que $U(L_0 - \phi(\lambda)) = p + V E(u)$.
Dans le cas (i), on remarque que $U$ et $p$ sont associés modulo $\Fil^1
S$ (\emph{i.e.} ils définissent le même idéal principal de $S/\Fil^1
S$), alors que dans le cas (ii), on observe qu'étant donné que les
termes constants de $U$ et $L_0 - \phi(\lambda)$ sont tous deux
divisibles par $p$, le terme constant de $V$ est, lui, congru à
$-c_0^{-1}$ modulo $p$. Finalement, remarquons que si $t$ satisfait les
conditions du lemme \ref{lem:petitt}, il en est de même de $\tr_1(t)$ et
qu'en outre modifier $t$ en $\tr_1(t)$ ne change pas la valeur de $v$.
Ainsi, on peut supposer sans perte de généralité que $t = \tr_1(t)$ ; 
c'est ce que nous faisons à partir de maintenant.

Soit $\calA \subset S$ l'idéal formé des éléments $A$ qui satisfont le
prémisse de $\H(\calL)$, c'est-à-dire :
$$\calA = \{ A \in S \, / \, A(Z - \calL_2) \in pS + \Fil^2 S_{K_0}\}.$$

\begin{lemme} 
\label{lem:genA}
L'idéal $\calA$ est engendré par $p + tE(u)$, $U E(u)$ et $\Fil^2 S$.
\end{lemme}

\begin{proof}
Bien entendu, on a $\Fil^2 S \subset \calA$. En outre, la formule
\eqref{eq:calA} combinée à la première condition du lemme \ref{lem:petitt}
montre que $p + t E(u) \in \calA$.

Soit $A \in \calA$ un élément arbitraire. Quitte à ajouter à $A$ un
élément de $\Fil^2 S$, on peut supposer que $A = \tr_2(A)$. Par le lemme
\ref{lem:A0}, on peut alors écrire $A = pA_0 + A_1 E(u)$ avec $A_0, A_1
\in W[u]$. Quitte à soustraire maintenant le multiple adéquat de $p +
tE(u)$, on peut supposer que $A_0 = 0$, et donc que $A = A_1 E(u)$. On
est finalement ramené à déterminer les polynômes $A_1 \in W[u]$ pour
lesquels $A_1 E(u) \in \calA$. Par une nouvelle application de la
formule \eqref{eq:calA}, c'est le cas si et seulement si $A_1 (L_0 -
\phi(\lambda)) \in pS + \Fil^1 S$. 

Dans le cas (i), $L_0 - \phi(\lambda)$ est inversible, d'où on obtient
l'équivalence entre $A_1 (L_0 - \phi(\lambda)) \in pS + \Fil^1 S$ et
$A_1 E(u) \in pE(u) S + \Fil^2 S$. Comme $U E(u)$ et $p E(u)$ sont
associés modulo $\Fil^2 S$, la condition est encore équivalente à $A_1
E(u) \in U E(u) S + \Fil^2 S$, et le lemme est démontré. 

Dans le cas (ii), on note que $p$ divise $A_1(\pi)(L_0(\pi) -
\phi(\lambda))$. On en déduit qu'il existe $U' \in S$ tel que $A_1 - U'
U \in \Fil^1 S$, puis que $A_1 E(u) \in U E(u) S + \Fil^2 S$ comme
annoncé. Réciproquement $U (L_0 - \phi(\lambda)) \in pS + \Fil^1 S$ par
construction de $U$, d'où $U E(u) \in \calA$.
\end{proof}

Pour tout $A \in \calA$, définissons
$$ m(A) = A \cdot f_1 + \frac1p \tr_2(A(\calL_2 - Z)) \cdot f_2 \,.$$
À partir de la congruence \eqref{eq:condfilr} et du lemme \ref{lem:genA},
on démontre rapidement le corollaire suivant.

\begin{cor}
Pour tout $A \in \calA$, $m(A)$ est un élément de $\Fil^2 \calM$. De
plus, $\Fil^2 \calM$ est engendré par $m(p+tE(u))$, $m(UE(u))$ et
$(\Fil^2 S)\cal M$.
\end{cor}

Dans la suite, nous écrirons $\overline{m}(A)$ à la place de
$\overline{m(A)}$. On cherche à présent à construire \emph{deux}
éléments qui engendrent $\Fil^2 \barM$ et à calculer leur image par
$\phi_2$. Le lemme suivant nous sera utile.

\begin{lemme}
\label{lem:phimodp}
Soit $A \in \calA$. Alors 
$$\phi_2(\overline{m}(A)) = \overline{(\phi(A)/p) f_1} + \overline{C
f_2}$$
où
$$C = \frac{\phi \circ \tr_2(A \calL_2) - \phi(A)Z +
\phi(\phi(\lambda))\phi(A - \tr_2(A))}{p^2} \in S.$$
\end{lemme}

\begin{proof} 
Le lemme suivra de la preuve de la proposition \ref{prop:condsuf} une
fois que l'on aura montré que $C$ et $\frac{1}{p^2}(p\phi(B) + \phi(AZ)
- \phi(A)Z)$ sont congrus modulo $p$, où $B = \frac{1}{p} \tr_2
(A(\calL_2- Z))$.

Puisque $A \in pS + \Fil^1S$ et $Z - \phi(\lambda) \in pS + \Fil^2S$, on
a $A(Z-\phi(\lambda)) \in pS + \Fil^3 S$. Ainsi
$$\tr_2(A(Z-\phi(\lambda))) - A(Z-\phi(\lambda)) \in p\Fil^2 S + \Fil^3 S$$
puis, en appliquant $\phi$ :
$$ \phi \circ \tr_2(A(Z-\phi(\lambda))) - \phi(A(Z-\phi(\lambda)))
\equiv 0 \pmod{p^3}.$$
Après un petit calcul :
$$p^2 C \equiv p\phi(B) + \phi(AZ) - \phi(A)Z \pmod{p^3}$$
ce qui est exactement ce que l'on voulait.
\end{proof}

\begin{prop}
\label{prop:genfil2}
\begin{enumerate}
\item Dans le cas (i), $\Fil^2 \barM$ est engendré par $\overline{m}
(p+tE(u))$ et $\overline{E(u)^2 f_1}$.

\item Dans le cas (ii), $\Fil^2 \barM$ est engendré par $\overline{m}
(p+tE(u))$ et $\overline{m}(U E(u))$.
\end{enumerate}
\end{prop}

\begin{proof}
Il suffit, dans chacun des cas, de montrer que les images par $\phi_2$
des deux éléments qui apparaissent dans l'énoncé du lemme engendrent
$\barM$ comme $\widetilde{S}_1$-module.

On commence avec $A = p + t E(u)$. Alors $\phi(A)/p = 1 + c\phi(t)$ et
$\overline{(\phi(A)/p)}$ est une unité in $\widetilde{S}_1$ grâce à la
seconde condition du lemme \ref{lem:petitt}. On calcule maintenant
l'élément $C$ du lemme \ref{lem:phimodp}. Comme $A \calL_2 \equiv pL_0 +
(tL_0 + L_1)E \pmod{\Fil^2 S_{K_0}}$, le premier terme du numérateur de
$C$ est
$$\phi \circ \tr_2 (A\calL_2) = p \phi(L_0) + pc \cdot \phi \circ
\tr_1(t L_0 + L_1)\,.$$
Le second terme, quant à lui, est $\phi(A)Z = p\phi(L_0) + pc \phi(t
\phi(\lambda))$.
Comme $\tr_2(A) = A$, le troisième terme s'annule et, après s'être 
rappelé que l'on a supposé $t = \tr_1(t)$, on obtient 
$$p^2 C = pc \cdot \phi\circ \tr_1( t L_0 + L_1 - t\phi(\lambda))\,.$$
La première condition du lemme \ref{lem:petitt} montre que la quantité
précédente s'annule, et donc que $C = 0$. Au final
$$\phi_2(\overline{m}(p+tE(u))) = \overline{(1 + c\phi(t)) f_1}.$$
On retiendra en particulier que le coefficient devant $\overline{f_1}$ 
dans le membre de droite est inversible.

Intéressons-nous maintenant à l'autre générateur, c'est-à-dire $m =
\overline{E(u)^2 f_1}$ dans le cas (i) et $m = \overline{m}(U E(u))$
dans le cas (ii). Dans le cas (i), on remarque que puisque $\tr_2
(E(u)^2(\calL_2 - Z)) =0$, on a $m = \overline m(A)$ avec $A = E(u)^2$.
Ainsi on peut appliquer le lemme \ref{lem:phimodp} pour calculer
$\phi_2(f)$. Comme $\phi(A) = p^2 c^2$, on a $\overline{\phi(A)/p} = 0$.
Par ailleurs, après avoir remarqué que $\phi\circ \tr_2 (A\calL_2) = 0$,
on trouve facilement
$$C = \frac{c^2 \phi(\phi(\lambda) - L_0)}{1 + c\phi(t)}.$$
Ainsi $\phi_2(m) = \overline{C f_2}$ et le fait qu'il engendre tout
$\barM$ avec $\phi_2(\overline{m}(p+tE(u)))$ résulte de l'inversibilité
de $C$ satisfaite puisque $\phi(\lambda) - L_0$ est une unité dans ce
cas.

Dans le cas (ii), on a clairement $m = \overline m(A)$ avec $A = U
E(u)$. On peut donc appliquer à nouveau le lemme \ref{lem:phimodp}. On a
d'abord clairement $\overline{\phi(A)/p} = \overline{c \phi(U)}$. Par
ailleurs, après un assez long calcul qui utilise les égalités $\tr_2(A)
= A$, $\tr_1(\phi(\lambda) U) = \phi(\lambda) U$ et $\tr_2(A \calL_2)
= \tr_1(U L_0) E(u) = (p + \phi(\lambda) U) E(u)$ ainsi que la
définition de $U$ et $V$, on trouve
$$ C = \frac{c^2 \phi(t - V)}{1 + c\phi(t)}\,.$$
Donc
$$\phi_2(\overline m(U E(u))) = \overline{c \phi(U) f_1} + \overline
{C f_2}.$$
Pour conclure, il suffit donc encore une fois de justifier que $C$ est
inversible dans $S$. Or nous avons déjà évalué les termes constants de
$t$ et $V$, et trouvé qu'ils sont respectivement congrus à
$-\frac{j}{ec_0}$ et $-\frac{1}{c_0}$ modulo $p$. Ainsi, ils ne sont
pas congrus entre eux, et la proposition s'ensuit.
\end{proof}

Définissons $g_1 = \overline{m}(p + tE(u))$ et $g_2 = \overline{E(u)^2
f_1}$ (resp. $g_2 = \overline{m}(UE(u))$) dans le cas (i) (resp. le cas 
(ii)). Posons également
$$B_1 = (L_0 - \phi(\lambda)) + \tr_1
\left(\frac{t(\phi(\lambda)-Z)}{p}\right) u^e $$
et dans le cas (ii) seulement
$$B_2 = \left(1 + \frac{1}{p} \tr_1 (U(\phi(\lambda)-Z))\right) u^e
\,.$$
On remarque tout de suite que par le lemme \ref{lem:Z}, $B_1$ et $B_2$
sont tous les deux des éléments de $S$. La structure de $\barM$ est
alors résumée dans la proposition suivante.

\begin{prop}
\begin{enumerate}
\item Dans le cas (i), $\Fil^2 \barM$ est engendré par
\begin{eqnarray*}
g_1 & =  & u^e \overline{t f_1} + \overline{B_1 f_2} \\
g_2 & =  & u^{2e} \overline{f_1}
\end{eqnarray*}
avec
\begin{eqnarray*}
\phi_2(g_1) & = & \overline{(1 + c\phi(t)) f_1} \\
\phi_2(g_2) & = & \overline{\frac{c^2\phi(\phi(\lambda)-L_0)}{1 +
c\phi(t)} f_2}
\end{eqnarray*}
où les deux coefficients dans les deux dernières équations sont des
unités.

\item Dans le cas (ii), $\Fil^2 \barM$ est engendré par
\begin{eqnarray*}
g_1 & =  & u^e \overline{t f_1} + \overline{B_1 f_2} \\
g_2 & =  & u^e \overline{U f_1} + \overline{B_2 f_2}
\end{eqnarray*}
avec
\begin{eqnarray*}
\phi_2(g_1) & = & \overline{(1 + c\phi(t)) f_1} \\
\phi_2(g_2) & = & \overline{c \phi(U) f_1} +
\overline{\frac{c^2\phi(t-V)}{1 + c\phi(t)} f_2}
\end{eqnarray*}
et les coefficients de $f_1$ dans $\phi_2(g_1)$ et $f_2$ dans
$\phi_2(g_2)$ sont des unités.
\end{enumerate}
\end{prop}

\begin{proof}
À la lumière de la proposition \ref{prop:genfil2} et de sa preuve, il ne
reste qu'à montrer que nos formules pour $B_1$ et $B_2$ sont correctes.

Pour le premier, il s'agit de montrer que si $A = p + tE(u)$ et $B =
\frac{1}{p} \tr_2(A(\calL_2 - Z))$, alors $B \equiv B_1 \pmod p$. Or, le
lemme \ref{lem:Z} donne $pZ \equiv p \phi(\lambda) \pmod{p^2S + \Fil^2
S}$, à partir de quoi un calcul direct utilisant la première condition
du lemme \ref{lem:petitt} entraîne
$$A(\calL_2 - Z) \equiv p(L_0 - \phi(\lambda)) + t E(u) (\phi(\lambda)
- Z) \pmod{p^2S + \Fil^2 S_{K_0}}$$
puis la conclusion.

Pour le second, il s'agit de même de montrer que si $A = U E(u)$ et
$B = \frac{1}{p} \tr_2(A(\calL_2 - Z))$, alors $B \equiv B_2 \pmod p$.
Mais c'est immédiat par le lemme \ref{lem:Z} et la définition de $U$.
\end{proof}

Dans le cas (i), on observe que $\overline{B_1}$ est inversible, alors
que $\overline{t}$ est soit égal à $0$, soit le produit de $u^j$ ($j \leq
1$) par une unité avec $v = j/e$ si $0 < j < e$ et $v \geq 1$ sinon. 
Dans le cas (ii), par contre, $\overline{B_1}$ (resp. $\overline{U}$,
resp. $\overline{B_2}$) est le produit d'une unité par $u^j$ (resp.
$u^{e-j}$, resp. $u^e$), alors que $\overline{t}$ est, quant à lui,
inversible (c'est-à-dire $v = 0$). De plus, on vérifie que le
coefficient en $u^e$ dans $\overline {tB_2} - \overline{B_1U}$ n'est
autre que le terme constant de $\overline{t}-\overline{V}$, dont on sait
qu'il ne s'annule pas. La première partie du théorème \ref{theo:red}
résulte facilement des considérations précédentes, alors que la seconde
est une conséquence du théorème 5.2.2 de \cite{caruso-crelle} et de la
proposition suivante.

\begin{prop}
Soit $\barM$ un objet de $\ModphiN{\widetilde{S}_1}$ de rang $2$
admettant pour base $(e_1, e_2)$. Dans la suite de la proposition,
toutes les lettres grecques font références à des éléments inversibles
de $\widetilde{S}_1$.

\begin{enumerate}
\item[0.] Supposons que $\Fil^2 \barM$ soit engendré par $\e_2$ et $u^{2e}
\e_1$ et que 
$$ \phi_2(\e_2) = \mu \e_1\,,\qquad \phi_2(u^{2e} \e_1) = \rho \e_2 \,.$$
Alors $T_{\st}(\barM)$ est irréductible et les pentes de son polygone
de l'inertie modérée sont $0$ et $2$.

\item Supposons que $\Fil^2 \barM$ soit engendré par 
$\alpha u^{e+j} \e_1 + \e_2$ et $u^{2e} \e_1$ et que 
$$ \phi_2(\alpha u^{e+j} \e_1 + \e_2) = \mu \e_1\,,\qquad \phi_2(u^{2e}
\e_1) = \rho \e_2 \,.$$
Si $j \geq e$, alors $T_{\st}(\barM)$ est irréductible et les pentes de
son polygone de l'inertie modérée sont $0$ et $2$. Au contraire, si $j <
e$, alors $\barM$ admet un sous-objet $\barM'$ de rang $1$ tel que
$\Fil^2 \barM' = u^{e-j} \barM'$.

\item Supposons que $\Fil^2 \barM$ soit engendré par $\alpha u^{e} \e_1
+ \beta u^j \e_2$ et $\gamma u^{2e-j} \e_1 + \delta u^e \e_2$ avec $j <
e$ et $\alpha \delta - \beta \gamma \in \widetilde{S}_1^{\times}$. 
Supposons également que
$$ \phi_2(\alpha u^{e} \e_1 + \beta u^j \e_2) = \mu \e_1\,,\qquad
\phi_2(\gamma u^{2e-j} \e_1 + \delta u^e \e_2) = \sigma u^{p(e-j)} \e_1
+ \rho \e_2 \,.$$
Alors $\barM$ admet un sous-objet $\barM'$ de rang $1$ tel que $\Fil^2
\barM' = u^e \barM'$.
\end{enumerate}
\end{prop}

\begin{proof}
L'alinéa (0) est immédiat par le théorème 5.2.2 de \cite{caruso-crelle}.  

Passons à (1). Si $j \geq e$, posons $\e_1' = \phi_2(\e_2) = \mu \e_1 -
\phi(\alpha) u^{p(j-e)} \e_2$. À partir de $e_2 \in \Fil^2 \barM$ et
$\phi(u^{2e}) = 0$, on obtient $\phi(u^{2e} \e_1') = \phi(\mu) \rho
e_2$. Comme $(\e_1', \e_2)$ est encore une base de $\barM$, le résultat
suit de (0).

Supposons désormais $j < e$.  On cherche $\barM'$ engendré par un
vecteur $m \in \barM$ de la forme $m = \e_2 + X u^{p(e-j)} \e_1$ où $X$
est inversible dans $\widetilde{S}_1$. On calcule
$$\phi_2(u^{e-j} \e_2) = -\rho \phi(\alpha) \e_2 + \mu u^{p(e-j)} \e_1 \,.$$
puis, en utilisant $2e < (p+1)(e-j)$ :
$$\phi_2(u^{e-j} m) = \rho (- \phi(\alpha) + \phi(X)
u^{p((p+1)(e-j)-2e)}) \e_2 + \mu u^{p(e-j)} \e_1.$$
Ainsi, si $X$ est une solution
\begin{equation}
\label{eq:X}
\rho X (- \phi(\alpha) + \phi(X) u^{p((p+1)(e-j)-2e)}) = \mu
\end{equation}
on a terminé. Mais on montre facilement que \eqref{eq:X} admet une
unique solution en résolvant l'équation coefficient par coefficient. En
outre, le coefficient constant de $X$ est celui de $-\frac \mu{\rho
\phi(\alpha)}$, d'où on déduit l'inversibilité souhaitée de $X$. Il 
faut encore vérifier que $\barM'$ est stable par $N$, mais cela ne 
pose pas de problème car la relation de commutation à $\phi_2$ 
implique directement $N\circ \phi_2(u^{e-j} m) = 0$ puis la stabilité
souhaitée étant donné que $\phi_2(u^{e-j} m)$ est un générateur de
$\barM'$.

Finalement, on traite (2). On cherche à nouveau $\barM'$ engendré par un
vecteur $m \in \barM$ de la forme $m = \e_2 + X u^{p(e-j)} \e_1$,
toutefois sans imposer à $X$ cette fois-ci d'être inversible.  En posant
$\Delta = \alpha \delta - \beta\gamma$, on calcule
$$\phi_2(\Delta u^e \e_2) = \rho \phi(\alpha) e_2 + (\phi(\alpha)\sigma
- \phi(\gamma)\mu) u^{p(e-j)} \e_1\,.$$
On a $u^{p(e-j)} \e_1 \in \Fil^2\barM$ car $p(e-j) > 2e$, d'où il suit
$\phi_2(u^e \cdot u^{p(e-j)}\e_1) = 0$. Ainsi
$$ \phi_2(u^e m) = \phi_2(u^e \e_2) = \frac{\rho
\phi(\alpha)}{\phi(\Delta)} \left(\e_2 + \frac{\phi(\alpha)\sigma -
\phi(\gamma)\mu}{\rho \phi(\alpha)} u^{p(e-j)} \e_1 \right)$$
et, puisque $\frac{\rho \phi(\alpha)}{\phi(\Delta)}$ est inversible, on
peut prendre $X = \frac{\phi(\alpha)\sigma - \phi(\gamma)\mu}{\rho
\phi(\alpha)}$. De même que dans le cas précédent, on vérifie pour finir
que $\barM'$ est stable par $N$.
\end{proof}

\end{document}